\documentclass[a4paper,fleqn]{cas-sc}

\usepackage{hyperref}
\usepackage{xcolor}
\usepackage{lipsum}
\usepackage{amsfonts}
\usepackage{graphicx}
\usepackage{epstopdf}
\usepackage{hhline}

\usepackage{multirow}%
\usepackage{amsmath,amssymb,amsfonts}%
\usepackage{amsthm}%
\usepackage{mathrsfs}%
\usepackage[title]{appendix}%
\usepackage{xcolor}%
\usepackage{textcomp}%
\usepackage{manyfoot}%
\usepackage{booktabs}%
\usepackage{algorithm}%
\usepackage{algorithmicx}%
\usepackage{algpseudocode}%
\usepackage{listings}%

\usepackage{booktabs}

\usepackage{booktabs}  
\usepackage{multicol}
\usepackage{bm}
\usepackage{algorithm}
\usepackage{rotating}
\usepackage{fancyvrb}
\usepackage{caption}

\usepackage{wrapfig}
\usepackage{subcaption}
\usepackage{color}
\usepackage{amsmath}
\usepackage{amssymb}
\usepackage{array}

\newtheorem{definition}{Definition}
\newtheorem{theorem}{Theorem}
\newtheorem{lemma}{Lemma}

\newtheorem{prop}{Proposition}

\usepackage{amsmath,amsfonts}

\usepackage[numbers,sort&compress]{natbib}

\usepackage{pifont}
\newcommand{\cmark}{\ding{51}}%
\newcommand{\xmark}{\ding{55}}%

\def\x{\mathbf{x}}
\def\c{\mathbf{c}}
\def\y{\mathbf{y}}

\begin{document}
\let\WriteBookmarks\relax
\def\floatpagepagefraction{1}
\def\textpagefraction{.001}

\shorttitle{}    

\shortauthors{}  

\title [mode = title]{An Efficient Algorithm for Vertex Enumeration of Arrangement}  

\author[1]{Zelin Dong}[]
\fnmark[1]
\ead{1155173731@link.cuhk.edu.hk}
\affiliation[1]{organization={Department of Mathematics, The Chinese University of Hong Kong, Hong Kong}}

\author[2]{Fenglei Fan}[]
\fnmark[2]
\ead{flfan@math.cuhk.edu.hk}
\affiliation[2]{organization={Department of Mathematics, The Chinese University of Hong Kong, Hong Kong}}

\author[3]{Huan Xiong}[]
\fnmark[3]
\ead{huan.xiong.math@gmail.com}
\affiliation[3]{organization={IASM, Harbin Institute of Technology, China}}

\author[4]{Tieyong Zeng}[]
\cormark[1]
\fnmark[4]
\ead{zeng@math.cuhk.edu.hk}
\affiliation[4]{organization={Department of Mathematics, The Chinese University of Hong Kong, Hong Kong}}

\cortext[1]{Corresponding author}

\begin{abstract}
This paper presents a state-of-the-art algorithm for the vertex enumeration problem of arrangements. We introduce a new pivot rule, called the Zero rule. The Zero rule possesses several desirable properties: i) It eliminates the objective function; ii) Its terminal dictionary is single; iii) We establish the if-and-only-if condition between the Zero pivot and its valid reverse pivot; iv) Applying the Zero pivot recursively definitely terminates in $d$ steps, where $d$ is the dimension of the input variables. {Based on this rule and its properties}, we leverage it to obtain a more efficient vertex enumeration algorithm. Theoretically, given an arrangement in $\mathbb{R}^d$ composed of $n$ hyperplanes with $v$ vertices, where $v_d$ of those vertices reach the terminal in exactly $d$ steps when applying the Zero pivot, the algorithm's complexity is $\mathcal{O}(n^2d^2(v-v_d)-ndv_d)$. {For the class of simple arrangements, the complexity is $\mathcal{O}(nd^4v)$, which significantly improves $\mathcal{O}(n^2dv)$ of the Avis and Fukuda algorithm}. Systematic and comprehensive experiments confirm that our algorithm is effective.
\end{abstract}


\maketitle

\section{Introduction}\label{sec1}

Vertex enumeration (VE) is a fundamental problem aimed at identifying vertices that satisfy specific linear constraints \cite{avis2010enumeration, bremner1997primal, dyer1983complexity, elbassioni2020enumerating, motzkin1953double, van1975node}. When these constraints take the form of linear inequalities, such as in the condition $Ax\leq b$, it involves enumerating vertices on a polytope, while the condition $Ax=b$ corresponds to enumerating vertices on an arrangement. Geometrically, VE bridges two crucial representations of a polytope or an arrangement: 
$H$-representation and $V$-representation, which respectively define a polytope or an arrangement from its faces and vertices. As such, VE is not only a fundamental problem in computational discrete geometry but also holds a wide array of real-world applications. VE has been applied in problems including, among others, effective field theory \cite{zhang2020convex} and computing approximating polytopes \cite{arya2022optimal}. More recently, VE has emerged as a key component in unraveling the inner workings of deep learning \cite{bergermann2021semi, fan2023deep, provan1994efficient, riihimaki2023simplicial, yang2021reachability, zago2023vertex}. The underlying concept stems from the fact that a network with piecewise linear activation functions serves as a piecewise linear entity, partitioning the space into numerous linear regions. Each region is defined by a set of linear constraints. Identifying the vertices of these linear constraints is pivotal in understanding the dynamics and robustness of deep networks \cite{zago2023vertex}.

The earliest published method to solve the problem of VE can be found in \cite{motzkin1953double}. Their approach is to add inequalities in a stepwise manner. In each step, "new" extreme points are created and "old" ones are excluded. Subsequently, due to the intrinsic connection between VE and linear programming, the simplex method and its variants, such as Avis and Fukuda's pivoting method \cite{avis1991pivoting}, Balinski's algorithm \cite{balinski1961algorithm}, Bremner's primal simplex method \cite{bremner1997primal}, Chand and Misra's algorithm \cite{chand1970algorithm}, Chvatal's linear programming approach \cite{chvatal1983linear}, Dyer and Frieze's complexity analysis \cite{dyer1983complexity}
, Mattheiss's approach \cite{Mattheiss1973AnAF},
Dyer, Martin E, and Proll's method \cite{dyer1977algorithm}
were utilized to solve the VE problems.
The simplex method was originally used in linear programming, based on the observation that a linear objective function always attains its maximum at extreme points of the feasible region (polytope). Geometrically, one only needs to move the current vertex to an adjacent one, and if the maximum exists, the algorithm will terminate at one vertex that maximizes the objective function. Specifically, in this method, a dictionary is used to represent the current vertex, and a step called "pivoting" is employed to select entering and leaving variables to update the dictionary, thus traversing between neighboring vertices. However, if a vertex lies on more than $d$ hyperplanes, the pivoting may trap into several dictionaries concerning the same vertex, which greatly deteriorates the efficiency. Hence, different pivot rules were introduced to avoid this situation and boost efficiency. For example, the Bland's rule \cite{bland1977combinatorial, bland1977new},  the Criss-Cross rule \cite{terlaky1985convergent,terlaky1987finite,wang1985finite} and its variant \cite{fukuda1991finiteness}, the Jensen's general relaxed recursion \cite{jensen1985coloring}, 
the Random-Facet pivot rule and its variants \cite{hansen2015improved,Kalai1992ASR, matouvsek1992subexponential}, the Random-Edge pivot rule \cite{hansen2014improved,matouvsek2006random}, and the Facet pivot rule \cite{yang2021facet,yang2022facet}. 

Many previous papers were dedicated to the VE of a polytope. Although the importance of VE of an arrangement is increasingly recognized, to the best of our knowledge, only a few methods focus on it: The method in Moss \cite{moss2012basis} is designed for listing exactly one co-basis from each orbit under the action from a set of automorphisms to the arrangement. It considers the smallest positive value and the largest negative value in the ratio test of every column, thereby obtaining all adjacent dictionaries, and hence all vertices. By repeating this process, every dictionary is obtained. Then, the algorithm compares each dictionary in this arrangement to determine whether it is a duplicate. This method is essentially a brute-force enumeration, with a high demand for storage and computational complexity. The method proposed by Avis and Fukuda \cite{avis1991pivoting} is based on an observation that since repeating one pivot rule can yield a path from some dictionary to an optimal dictionary, then tracking back the path can generate all dictionaries ending at that optimal dictionary, without the need for the brute-force enumeration. Instead of directly computing coordinates of all vertices, the algorithm first outputs lexicographically minimum basis for each vertex, and then computes coordinates. 

Despite the attraction of the idea, the algorithm in Avis and Fukuda \cite{avis1991pivoting} could be upgraded from the following three major aspects: i) The method for examining valid reverse pivots takes several steps, which is less efficient than an explicit if-and-only-if condition according to our analysis in Section \ref{sec:ca} and \hyperref[analysis_iff]{Appendix}. ii) The objective function is not necessary for the vertex enumeration of an arrangement. Firstly, the objective function restricts the traversal of vertices, resulting in more steps in complex arrangements. Secondly, the value of the objective function is unpredictable until a vertex is reached, which means that the objective function cannot be optimized until visiting all vertices. This is more time-consuming. Thirdly, the pivot operation on a dictionary can be performed without the objective function, and neither the magnitude of the objective function nor the function itself affects computing the vertex coordinates based on the dictionary's basis or co-basis. Therefore, the objective function merely guides the pivot process, however, this guidance can be directly captured from the dictionary itself. iii) Multiple optimal vertices or dictionaries may exist in the Criss-Cross rule, which introduces additional steps and storage.

In this paper, we first introduce a pivot rule that can be applied throughout the hyperplane arrangement, referred to as the Zero rule, which utilizes the linear independence of the normal vectors of hyperplanes in the arrangement. Compared to the Criss-Cross rule, its most notable features are: i) {the number of Zero pivots from each dictionary to the terminal is upper bounded by the dimension $d$, which is a small constant bound};
ii) it eliminates the objective function; iii) the collection of dictionaries from which the Zero rule cannot select a suitable position is always unique, while the Criss-Cross rule may have multiple such dictionaries.
{In addition, we also give an if-and-only-if condition for detecting if a pivot is a valid reverse Zero pivot.}


Based on these desirable properties of the Zero rule, we realize that the key points of the algorithm proposed by Avis and Fukuda \cite{avis1991pivoting} can be further improved with the help of the Zero rule. Compared to the algorithm of Avis and Fukuda \cite{avis1991pivoting}, using our algorithm can effectively reduce both space and time complexity. {Specifically, given a hyperplane arrangement in $\mathbb{R}^d$ consisting of $n$ distinct hyperplanes with $v$ vertices, there are $\mathcal{O}(v)$ dictionaries in this arrangement. Note that applying the Zero pivot takes at most $d$ iterations. Therefore, we denote the number of dictionaries that require exactly $d$ iterations as $\mathcal{O}(v_d)$. Then the computational complexity of the algorithm improved by the Zero rule is $O(n^2 d^2 (v - v_d) + nd v_d)$. If a hyperplane arrangement is a simple arrangement \cite{toth2004handbook}, where $v=\binom{n}{d}$ and $v_d \geq\binom{n-d}{d}$, the complexity favorably turns into at most $\mathcal{O}(nd^4v)$, which significantly improves Moss \cite{moss2012basis} and Avis \& Fukuda \cite{avis1991pivoting}. Moreover, since there is no additional storage, the needed
space are for those variables generated during the computation, which are at most two dictionaries, in the order of $\mathcal{O}(nd)$. Table \ref{tab:summary} summarizes the information of the proposed algorithm and its counterparts.}

\begin{table*}[htbp]
\renewcommand{\arraystretch}{1.5}
  \centering
  \scalebox{0.9}{
  \begin{tabular}{|c|c|c|c|c|c|}
    \hline & storage & obj. func.  & complexity & single tree\\ \hline 
    Ours  & $\mathcal{O}(nd)$ & \xmark  & $\mathcal{O}(n^2 d^2 (v - v_d) + nd v_d)$ & \cmark \\ \hline
    \texttt{AF} \cite{avis1991pivoting} &  $\mathcal{O}(nd)$ & \cmark & $\mathcal{O}(n^2dv)$ & \xmark \\ \hline
    \texttt{Enhanced AF} &  $\mathcal{O}(nd)$ & \cmark & $\mathcal{O}(n^2dv)$ & \xmark \\ \hline  
    \texttt{Moss} \cite{moss2012basis} & $\mathcal{O}(ndv)$ & \xmark & $\mathcal{O}(nd^2v^2)$ & \cmark\\ \hline
  \end{tabular}}
  \label{results:sketch}
  \caption{{A summary between the proposed algorithm, \texttt{Enhanced AF}, Avis and Fukuda's algorithm (\texttt{AF} \cite{avis1991pivoting}), Moss' algorithm \cite{moss2012basis} (\texttt{Moss}). In our algorithm, $v_d$ is the number of dictionaries that require exactly $d$ steps of pivoting. If a hyperplane arrangement is a simple arrangement, the complexity of our algorithm is $\mathcal{O}(nd^4v)$, while for a hyperplane arrangement where $v_d$ accounts for a small portion, the complexity is $\mathcal{O}(n^2d^2v)$. Fortunately, the simple arrangement is quite general.
  }}
  \label{tab:summary}
  \vspace{-0.3cm}
\end{table*}

{The structure of this paper is as follows: Section \ref{sec2} introduces the definitions and notations used throughout this paper. Section \ref{sec3} elaborates on the important properties of the proposed Zero rule. Section \ref{sec4} includes: 1) a brief description of the algorithm proposed by Avis and Fukuda \cite{avis1991pivoting}; 2) 
the VE algorithm use the Zero Rule and advantage analysis, which is our algorithm; 3) a discussion of the complexities of various vertex enumerating algorithms. Section \ref{sec5} presents detailed toy examples. Section \ref{sec6} covers systematic experiments. Section \ref{sec7} concludes this paper. In addition to the main content, we also give an enhanced version to the \texttt{AF} algorithm for independent interest to some readers in the \hyperref[analysis_iff]{Appendix}, we denote the enhanced version as \texttt{Enhanced AF}.}

\section{Notations and Definitions}\label{sec2}

For convenience and consistency, notations in this paper are the same as \cite{avis1991pivoting} that uses dictionaries instead of the simplex tableau.

\begin{definition}[Hyperplane and arrangement]
    A hyperplane in $\mathbb{R}^d$ is defined as \\ $\{ y\in \mathbb{R}^d : \langle \c,\y \rangle=c_1y_1+ c_2y_2+ \ldots+ c_d y_d= b, ~ \c\in \mathbb{R}^d \}$.    
    The finite collection of hyperplanes is called hyperplane arrangement.
\end{definition}

\begin{definition}[Vertex of an arrangement]
    A vertex of an arrangement is the unique solution to the system of $d$ equations corresponding to $d$ intersecting hyperplanes. Moreover, a vertex is degenerate if it is contained in more than $d$ hyperplanes. A vertex is also called a basic solution in linear programming.
\end{definition}

In this paper, we assume that the arrangement contains at least one vertex. {Particularly, simple arrangements are those hyperplane arrangements such that every vertex is non-degenerate and every $d$ hyperplanes form a vertex. This definition can also be found on page 530 of the book \cite{toth2004handbook}.}

\begin{definition}[Slack variable]
    For each hyperplane, $\langle \c_i, \y \rangle = b_i$, the corresponding slack variable is defined as $x_i = b_i x_g - \langle \c_i, \y \rangle$, where $x_g$ is a constant-valued variable that is always equal to $1$.
\end{definition}

{It is important to note that if there are $d$ hyperplanes in the arrangement, $\langle\c_{j_i},y_0\rangle=b_{j_i}, ~i=1,\cdots, d$, such that their coefficients $\{\c_{j_1}, \cdots,\c_{j_{d}}\}$ are linearly independent, then by solving the system of equations $x_{j_i} = b_{j_i} x_g - \langle \c_{j_i}, \y \rangle$, we can express $\y$ in terms of $x_g$ and those slack variables. In this case, if we set $x_{j_i}=0, ~ \forall i=1,\ldots,d$, meaning that the corresponding $\y$ lies on all $d$ hyperplanes simultaneously, then the point $\y$ is a basic solution (vertex).}

\begin{definition}[Dictionary with an objective function]\label{def:dic_withobj}
{Using the above notations, assume there are $n$ hyperplanes in an arrangement. Let $x_f$ be the objective function, which is a linear combination of elements in $\y$. 
If $\{\c_{j_1}, \cdots,\c_{j_{d}}\}$ is linearly independent, express $y_1, \cdots, y_d$ in terms of $x_g, x_{j_1}, \cdots, x_{j_d}$. Consequently, put $\{i_1,\cdots,i_{n-d}\}=\{1,\cdots,n\}\backslash\{j_1,\cdots,j_d\}$, by plugging this representation into the objective function and other slack variables, we have
\begin{equation}
    \x_B = \bar{A}\x_N,
\label{eqn:dic}    
\end{equation}
where $\x_B = \{x_{i_1}, \cdots, x_{i_{n-d}}, x_f\}$, $\x_N = \{x_g, x_{j_1}, \cdots, x_{j_{d}}\}$. Such an equation system is called a dictionary. Sometimes we also denote a dictionary by $D$ or $(B, N, \bar{A})$.}
\end{definition}

\begin{definition}[Dictionary and Basic solution]
    {Using the above notations, setting $x_{j_1},\cdots,x_{j_d}$ to zeros induces a basic solution $\y_0$, we say $\y_0$ be corresponding basic solution of the dictionary $\x_B=\bar{A}\x_N$. In particular, if $\y_0$ is a degenerate vertex, we say $\x_B=\bar{A}\x_N$ a degenerated dictionary.}
\end{definition}

\begin{definition}[Basis and co-basis]
    In a dictionary $\x_B=\bar{A}\x_N$, $\x_B$ are basic variables, the collection of their indices, $B=\{i_1,\cdots,i_{n-d},f\}$ is called basis. $\x_N$ are non-basic variables, the collection of their indices, $N=\{g,j_1,j_2,\cdots, j_d\}$ is called co-basis. $\bar{A}=(\bar{a}_{ij})$ is the coefficient matrix, and $\bar{a}_{ij}$ expresses entries in $\bar{A}$ corresponding to $x_{i}\in \x_B$ and $x_{j}\in\x_N$.
\end{definition}

{In some case, we need to consider basic variables without the objective function or nonbasic variables without $x_g$. Thus, we put $B_{\neq f}=B\backslash\{f\}=\{i_1,\cdots,i_{n-d}\}$ and $N_{\neq g}=N\backslash\{g\}=\{j_1,j_2,\cdots, j_d\}$.}

\begin{definition}[Dictionary without the objective function]\label{def:dic_withoutobj}
By removing the row of $x_f$ in Eq. \eqref{eqn:dic}, we say it is a dictionary without objective function. In this case, $\x_B=(x_{i_1},\cdots,x_{i_{n-d}})$ and $B=B_{\neq f}=\{i_1,\cdots,i_{n-d}\}$.
\end{definition}

\begin{figure}[htbp]
    \centering
    \includegraphics[width=0.55\linewidth]{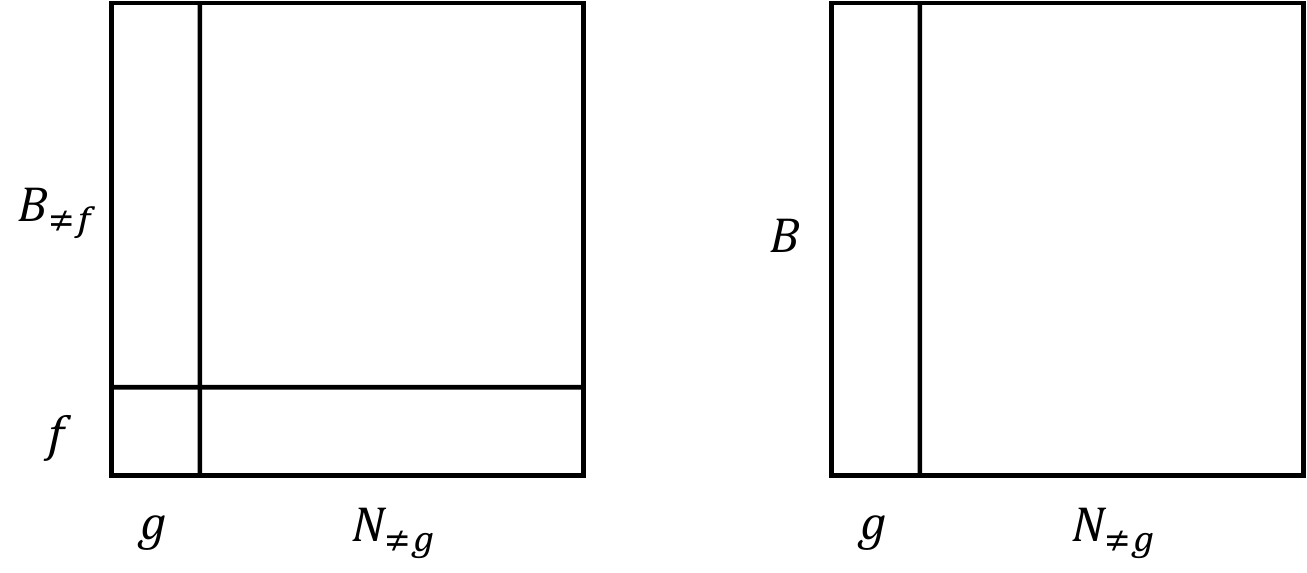}
    \caption{A dictionary with and without the objective function.}
    \label{fig:eg-dic}
    \vspace{-0.5cm}
\end{figure}

\begin{definition}[Primal Feasible]
Let $\x_B=\bar{A}\x_N$ be a dictionary. If $\bar{a}_{ig} \geq 0$, for some $i\in B_{\neq f} $, then the variable $x_i$ is primal feasible. A dictionary is \textit{primal feasible}, if $x_i$ is primal feasible $\forall i\in B_{\neq f}$,.
\end{definition} 

\begin{definition}[Dual Feasible and Optimal]
Let $\x_B=\bar{A}\x_N$ be a dictionary with the objective function. If $\bar{a}_{fj}\leq 0$ for some $j\in N_{\neq g}$, then the variable $x_j$ is dual feasible. A dictionary is \textit{dual feasible}, if $x_j$ is dual feasible $\forall j\in N_{\neq g}$. The dictionary is \textit{optimal}, if it is both primal and dual feasible.
\end{definition} 

\begin{definition}[Pivoting rule] {The pivoting rule is a method of selecting a pair of variables on a dictionary for subsequent operations.}
\end{definition}

\begin{definition}[The Criss-Cross rule]
Let $\x_B=\bar{A}\x_N$ be a dictionary with an objective function. The Criss-Cross rule \cite{terlaky1985convergent, terlaky1987finite, wang1985finite} selects the variables $(r,s)$ in the following order:

1) $i \neq f, g$ is the smallest index such that $x_i$ is primal or dual infeasible.

2) If $i \in B_{\neq f}$, let $r = i$ and $s$ be the minimum index such that $\bar{a}_{rs}> 0$; otherwise,
let $s = i$ and let $r$ be the minimum index such that $\bar{a}_{rs}< 0$.
\end{definition}

If the step $1)$ does not apply, then the dictionary is optimal.

\begin{definition}[Pivot]
    Let $\x_B=\bar{A}\x_N$ be a dictionary, a pivot $(r,s)$ on the dictionary is an exchange between $x_r$ in $\x_B$ and $x_s$ in $\x_N$, which generates a new dictionary $\x_{\tilde{B}}=\tilde{A}\x_{\tilde{N}}$, where the new coefficient matrix $\tilde{A}=(\tilde{a}_{ij})$ is obtained by computing Eq. \eqref{eqn:pivoting} and sorting rows and columns of  $\tilde{A}=(\tilde{a}_{ij})$ in increasing the order of indices:
    \begin{equation}
\begin{aligned}
    &\tilde{a}_{sr}=\cfrac{1}{\bar{a}_{rs}},&&\tilde{a}_{ir}=\cfrac{\bar{a}_{is}}{\bar{a}_{rs}},&&\tilde{a}_{sj}=-\cfrac{\bar{a}_{rj}}{\bar{a}_{rs}},&&\tilde{a}_{ij}=\bar{a}_{ij}-\cfrac{\bar{a}_{is}\bar{a}_{rj}}{\bar{a}_{sr}},&& i\in B_{\neq r},  j\in N_{\neq s}.
    \label{eqn:pivoting}
\end{aligned}
\end{equation}
In particular, if the pair $(r,s)$ on the dictionary $D$ is obtained by the $\mathcal{R}$ rule, we call the pivot $(r,s)$ an $\mathcal{R}$ pivot. 
\end{definition} 

{In fact, the pivot $(r,s)$ on a dictionary is moving $x_s$ to the left hand side and $x_r$ to the right hand side in $\x_B = \bar{A}\x_N$. This equation is to update the coefficient matrix and then rearrange the variables in a new basis and co-basis in lexicographic order to create a new dictionary. }

\begin{definition}[Valid Reverse Pivot]
    {Let $D$ be a dictionary obtained by a pivot $(r,s)$ on another dictionary $\tilde{D}$. Then the pivot $(s,r)$ on $D$ is called a valid reverse pivot. In particular, if the pivot $(r,s)$ on the dictionary $\tilde{D}$ is a $\mathcal{R}$ pivot, the pivot $(s,r)$ on $D$ is called a valid reverse $\mathcal{R}$ pivot.}
\end{definition}

According to the previous formulation, the vertex enumeration problem can be transformed into enumerating all dictionaries. The key in this process is to circumvent yielding vertices associated with multiple dictionaries more than once. Section \ref{sec4} gives the method that employs the reverse idea for vertex enumeration.

\section{The Zero rule}\label{sec3}

In this section, we propose the Zero rule which works on the dictionary $\x_B=\bar{A}\x_N$ without the objective function. Let $\x_B=\bar{A}\x_N$ be a dictionary without the objective function. Then the Zero rule is given in the following:

    1) Let $s$ be the smallest index in $N_{\neq g}$ such that $\{ i\in B ~|~ \bar{a}_{is}\neq0 \textup{ and } i<s \}\neq \emptyset$.

    2) After $s$ is selected, let $r=\min\{ i\in B ~|~ \bar{a}_{is}\neq0 \textup{ and } i<s \}$. 

\textbf{Remark 1:} The Zero rule refers to a selection process that does not impose any requirements for primal feasibility or dual feasibility. Consequently, the dictionary where the Zero rule can not select a proper position could be neither primal feasible nor dual feasible dictionary. For convenience, we refer to it as a \textit{terminal dictionary}, or simply as the \textit{terminal}. {The existence of the terminal dictionary is trivial. Note that any Zero pivot removes $s$ from the co-basis and introduces $r$ into it. Since $r<s$, after the Zero pivot, the new co-basis has a smaller lexicographical order. Therefore, the dictionary equipped with the smallest lexicographical order co-basis automatically becomes a terminal dictionary.}

Now, we prove the desirable properties of the proposed Zero rule. First, we show that given a dictionary, the Zero rule can only select a unique pair of variables to pivot, which is obtained by direct computation. 

\begin{prop}[Uniqueness]
    The Zero rule only selects a unique entry of a dictionary to pivot.
\end{prop}

\noindent\textit{Proof: }
Direct calculation.

\rightline{$\Box$}

\begin{figure}[htbp]
    \vspace{-0.5cm}
    \centering    \includegraphics[width=0.8\linewidth]{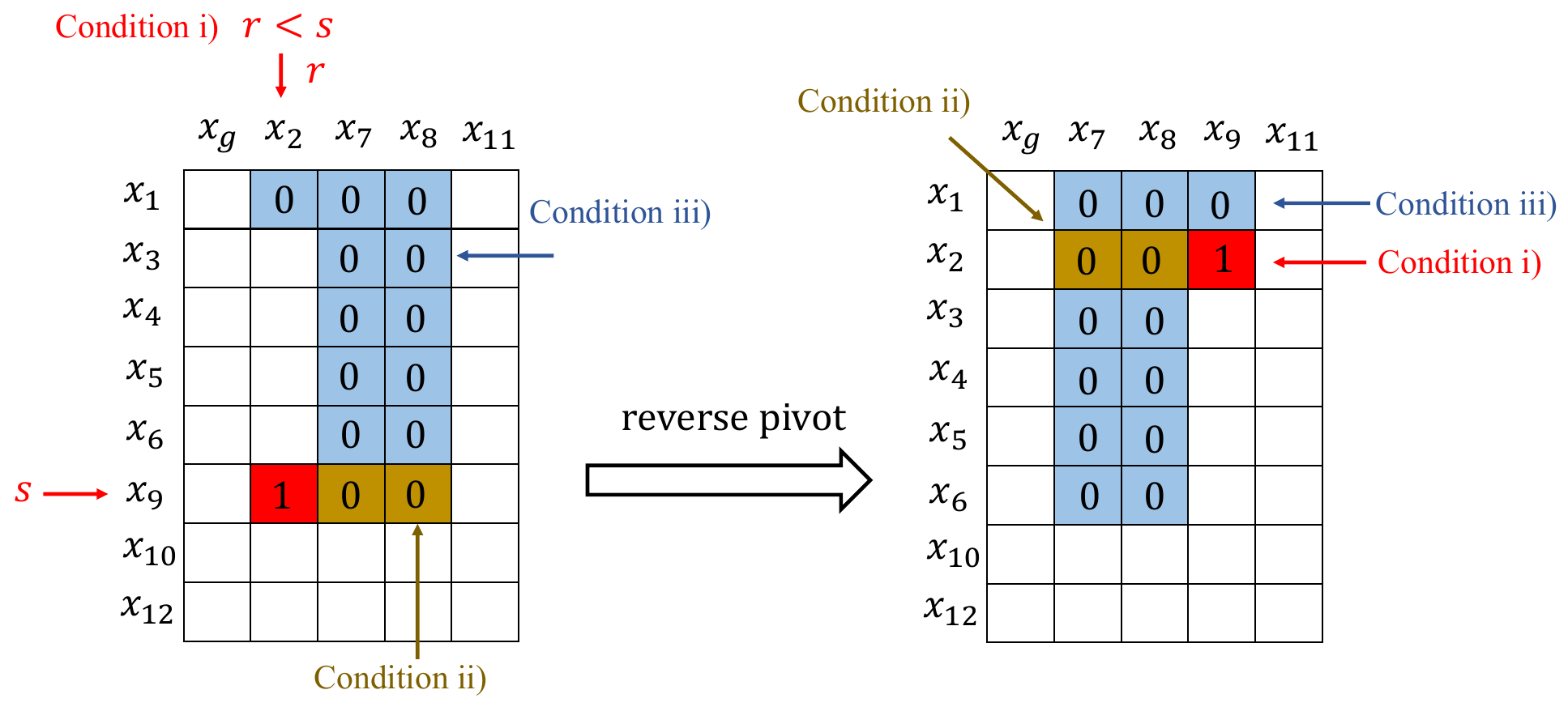}
    \caption{An exemplary illustration of the Zero rule and its valid reverse. Blocks of different colors represent where different sub-conditions are enforced.}
    \label{fig:enter-label}
    \vspace{-0.2cm}
\end{figure}

{Next, in order to have a better understanding of Zero pivots,} we give an if-and-only-if condition to test if a pivot on a dictionary is a valid reverse Zero pivot. Figure \ref{fig:enter-label} presents an example to illustrate the Zero rule and its valid reverse, where blocks of different colors correspond to different sub-conditions taking effect.



\begin{prop}[Reversibility] 
Let $\x_B=\bar{A}\x_N$ be an arbitrary dictionary. Then $(s,r)$, where $s\in B$ and $r\in N_{\neq g}$, is the valid reverse of the Zero rule if and only if 

i) $ r<s$ and $\bar{a}_{sr}\neq0$; 

ii) $\forall j\in N_{\neq g}$ and $r<j<s$, $\bar{a}_{sj}=0$;

iii) $\forall j\in N_{\neq g}$ and  $j<s$, for those $i\in B$
and $i<j$, $\bar{a}_{ij}=0$.
\label{prop:reverse}
\end{prop} 

\noindent\textit{Proof: }
Let $\x_{\tilde{B}}=\tilde{A}\x_{\tilde{N}}$ be the dictionary obtained from $\x_B=\bar{A}\x_N$ through pivoting $(s,r)$. Now, let us prove the if-and-only-if condition, respectively.

$\Rightarrow$: Assume $(s,r)$ is the valid reverse pivot of the proposed Zero rule, the resultant dictionary $\x_{\tilde{B}}=\tilde{A}\x_{\tilde{N}}$ will make the proposed Zero rule automatically pinpoint $(r,s)$ to pivot. This means that 1) $\forall j\in \tilde{N}_{\neq g}$ and $j<s$, any $i\in \tilde{B}_{\neq r}$ and $i<j$, $\tilde{a}_{ij}=0$; 2) $\tilde{a}_{rj}=0$ for any $r<j<s$; otherwise $(r,j)$ will be selected to pivot. 3) $\tilde{a}_{rs}\neq0$ and $r<s$. Now, let us use these three properties to deduce conditions i)-iii). 

{\tiny $\bullet$} The condition i) holds from 3) directly.

{\tiny $\bullet$} As far as the condition ii) is concerned, according to 2), we have $\tilde{a}_{rj}=0$ for any $r<j<s$. Combining $\tilde{a}_{rj}=\cfrac{\bar{a}_{sj}}{\bar{a}_{sr}}=0$ leads to that $\bar{a}_{sj}=0$ for $r<j<s$.

{\tiny $\bullet$} 1) is almost the condition iii). Combining Eq. \eqref{eqn:pivoting} and ii), $\tilde{a}_{ij}=0$ can naturally lead to $\bar{a}_{ij}=0$.

$\Leftarrow$: Assume conditions i)-iii) hold, we have

{\tiny $\bullet$} The condition i) ensures the entry in $(r,s)$ in $\x_{\tilde{B}}=\tilde{A}\x_{\tilde{N}}$ is non-zero and $r<s$.

{\tiny $\bullet$} As far as 3) is concerned, according to the condition ii), we have $\bar{a}_{sj}=0$ for any $r<j<s$. Combining $\tilde{a}_{rj}=\cfrac{\bar{a}_{sj}}{\bar{a}_{sr}}=0$ leads to that $\tilde{a}_{sj}=0$ for $r<j<s$.

{\tiny $\bullet$} The condition iii) is almost 1). By Eq. \eqref{eqn:pivoting} and ii), $\bar{a}_{ij}=0$ can naturally lead to $\tilde{a}_{ij}=0$.

\rightline{$\Box$}

After deriving the valid reverse of the Zero rule, we introduce another feature of indices leaving the co-basis when repetitively enforcing the Zero rule.

\begin{lemma}[Index increasing]
    Starting with any dictionary in any arrangement, when the Zero pivot is repeated, the sequence of indices leaving the co-basis in each step increases.
\label{lemma:index_increasing}    
\end{lemma}

\noindent\textit{Proof:} 
Assume in the $l$-th step, the Zero rule select $(r_l,s_l)$ and yields the $l$-th dictionary $\x_{B_l}=\bar{A}\x_{N_l}$, where the Zero rule pinpoints $(r_{l+1},s_{l+1})$.
$$\cdots\xrightarrow[Zero \ rule]{(r_l,s_l)} \x_{B_l}=\bar{A}_l\x_{N_l} \xrightarrow[Zero \ rule]{(r_{l+1},s_{l+1})} \cdots$$

Note that $(s_l,r_l)$ is the valid reverse pivot of the proposed Zero rule, then by Proposition \ref{prop:reverse}, we have $\forall j\in N_{\neq g} \textup{ and } j<s_l, \textup{ for those } i\in B
\textup{ and } i<j, \bar{a}_{ij}=0$.

For columns in $\x_{B_l}=\bar{A}\x_{N_l}$, from ii, iii), $\forall j\in N_{\neq g} \textup{ and } j<s_l$, $\{i\in B_l~|~i<j, \bar{a}_{ij}\neq0\}=\emptyset$. At the same time, $x_{s_l}$ has been moved into $\x_B$. Thus, applying the Zero rule again would cause the later column $s_{l+1}>s_l$ to be selected.
Hence, the index leaving the co-basis is increasing. 

\rightline{$\Box$}

Based on Lemma \ref{lemma:index_increasing}, one can prove that the number of Zero pivots recursively is no more than $d$. 

\begin{prop}[At most $d$ steps]
For any dictionary in any arrangement of $\mathbb{R}^d$, consecutively applying the Zero pivot always attain its terminal in at most $d$ steps.
\label{prop: d-steps conv}
\end{prop}

\noindent\textit{Proof: } According to Lemma \ref{lemma:index_increasing}, it can be concluded that each time the Zero rule is applied, the column it selects will be at least one position later compared to the previous step. In any dictionary in $\mathbb{R}^d$, there are a total of $d$ columns, excluding the column corresponding to $x_g$. This implies that the Zero rule will be repeated for at most of $d$ times. 

\rightline{$\Box$}

The property given in Proposition \ref{prop: d-steps conv} is very important because most pivot rules, such as Bland's rule, the Criss-Cross rule, the Random-Edge pivot rule, the Random-Facet pivot rule, converge in exponential or subexponential steps. The exponential convergence of Bland's rule is shown in \cite{Paparrizos2009LinearPK, terlaky1993pivot}, while the Criss-Cross rule is shown by an example that requires at least $2^n-1$ steps for convergence \cite{Roos1990AnEE}. Moreover, even better method, such as 
the Random-Edge pivot rule, has an exponential upper bound on pivot steps \cite{gartner2007two}, while the Random-Facet pivot rule in \cite{Hansen2015AnIV} has a subexponential upper bound. In comparison, our Zero rule exhibits a very fast convergence speed. {Additionally, for different hyperplane arrangements with the same $d$ and $n$, their upper bounds are different, making their pivot difficult to study. However, our Zero pivot has a common upper bound for each hyperplane arrangement in the same $d$, which is more convenient.}

After studying the number of Zero pivots, we focus on the number of terminals in the arrangement. Will there exist multiple terminal dictionaries? To address this question, it is necessary to identify what "0” means in each dictionary.

\begin{lemma} 
Let $\x_B=\bar{A}\x_N$ be an arbitrary dictionary. Then for any $i\in B$, the normal vector of the $i$-th hyperplane $\c_i\in \mathbb{R}^d$ satisfy $\c_i\in \mathrm{span}(\{ \c_j~|~\bar{a}_{ij}\neq0\})$.
\label{lemma: zeros}
\end{lemma} 

\noindent\textit{Proof}: In any fixed dictionary $x_B=\bar{A}x_N$, $x_i=\sum\limits_{j\in N}\bar{a}_{ij}x_{j}$. Also, note that we have 
\begin{equation}
    \begin{cases}
        & x_i=b_i x_g-\langle \c_i,\y \rangle, \\
        & x_j= b_j x_g-\langle \c_j,\y \rangle,
    \end{cases}
    \label{eqn:ij_relation}
\end{equation}
for all $i\in B $ and $j\in N_{\neq g}$. Substituting Eq. \eqref{eqn:ij_relation} into $x_i=\sum\limits_{j\in N}\bar{a}_{ij}x_{j}$, we can derive the relationship between $\c_i$ and $\c_j$ as follows:
\begin{equation}
    c_{il}=\sum\limits_{j\in N_{\neq g}}\bar{a}_{ij}c_{jl}, \qquad \qquad l\in \{1,2,...,d\},
\end{equation}
where $c_{il}$ and $c_{jl}$ are the $l$-th coordinate of the vector $\c_i$ and $\c_j$, respectively.  Wrapping it into a vector representation, we have

\begin{equation}
\begin{aligned}
    \c_i=\sum\limits_{j\in N_{\neq g}}(\bar{a}_{ij}\c_{j})=\sum\limits_{j\in N_{\neq g}, \bar{a}_{ij}\neq0\ } (\bar{a}_{ij}\c_{j}), 
    \label{eqn:lem}
\end{aligned}
\end{equation}
which means $\c_i\in \mathrm{span}(\{ \c_j~|~ \bar{a}_{ij}\neq0\})$.

\rightline{$\Box$}

From the above lemma, it can be seen that the proposed Zero rule is not an indexing game. Instead, it is deeply related to the properties of the arrangement itself, as represented by the dictionary. The relationship between zero and non-zero entries is the reason why we call our proposed rule the Zero rule. Now, we show that under the Zero rule, there only exists one terminal.

\begin{prop}[Unique terminal]
{The collection of terminal dictionaries of any hyperplane arrangement is single.}
\label{prop:convergence}
\end{prop}

\noindent\textit{Proof}: 
Let $\x_{B}=\bar{A} \x_{N}$ be the dictionary such that $N$ has attained the smallest lexicographical order among all co-basis in the hyperplane arrangement; note that the dictionary must be a terminal of the Zero rule. Furthermore, suppose there exists another terminal dictionary, namely  $\x_{\tilde{B}}=\tilde{A}\x_{\tilde{N}}$ and $\tilde{N}\neq N$, that is, 
\begin{equation}
    \begin{aligned}
    \forall s\in \tilde{N}_{\neq g}, \forall i<s \textup{ and } i\in \tilde{B}, \tilde{a}_{is}=0. \label{eqn:term}
    \end{aligned}
\end{equation}

Comparing $N$ and $\tilde{N}$ entry by entry, let the $j$-th pair be the first such that they are not equal, and denote them as $n_j\in N, \tilde{n}_j\in \tilde{N}$.

Then, it is possible to construct a partition on $N$ and $\tilde{N}$, given as $N=N_0\cup\{n_j\}\cup N_1, \tilde{N}=\tilde{N}_0\cup\{\tilde{n}_j\}\cup \tilde{N}_1$ where $N_0, \tilde{N}_0$ are the set of all numbers before $n_j$ and $\tilde{n}_j$ in $N$ and $\tilde{N}$ respectively. By the definition of $n_j$ and $\tilde{n}_j$, $N_0=\tilde{N}_0$. $N_1, \tilde{N}_1$ are remainders in $N$ and $\tilde{N}$, respectively.

\begin{figure}[htbp]
\vspace{-0.4cm}
    \centering
    \includegraphics[scale=0.4]{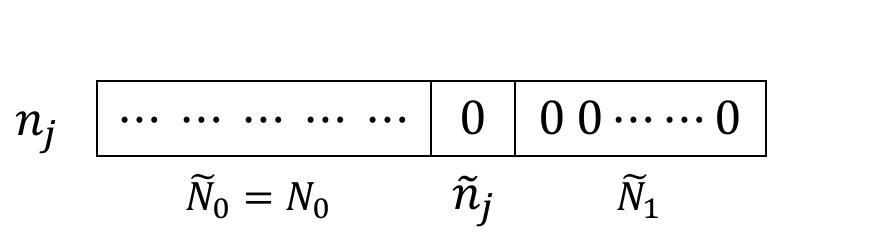}
    \caption{Zeros in the row of $x_{n_j}$ on $\x_{\tilde{B}}=\tilde{A}\x_{\tilde{N}}$.}
    \vspace{-0.2cm}
    \label{expalain:unique}
\end{figure}

Note that $\tilde{N}\neq N$, and $N$ attains the lexicographically minimum among all co-basis in the system, then $n_j<\tilde{n}_j$, and $\forall \tilde{n}_l \in \tilde{N}_1, n_j<\tilde{n}_j<\tilde{n}_l$. In the dictionary $\x_{\tilde{B}}=\tilde{A}\x_{\tilde{N}}$ by Eq. \eqref{eqn:term}, $\bar{a}_{n_{j}\tilde{n}_{j}}=0$ and for all $\tilde{n}_l \in \tilde{N}_1$, $\tilde{a}_{n_{j}\tilde{n}_{l}}=0$, which is shown as the Figure \ref{expalain:unique}.

It follows that $\{ \c_{\tilde{n}_{l}}~|~ \tilde{a}_{n_{j}\tilde{n}_{l}}\neq0\}\subseteq\{ \c_{\tilde{n}_{l}} ~|~ \tilde{n}_{l}\in N_0\}$. By Eq. \eqref{eqn:lem}, we have 
\begin{equation}
    \begin{aligned}
    \c_{n_{j}}\in \mathrm{span}(\{ \c_{\tilde{n}_{l}} ~|~ \bar{a}_{n_{j}\tilde{n}_{l}}\neq0\})\subseteq \mathrm{span}(\{ \c_{l} ~|~ l\in N_0\}).
    \label{eqn:relationship}
\end{aligned} 
\end{equation}

Since $N$ is a co-basis in the arrangement, then $\{ \c_{l} ~|~ l\in N_0\}\cup\{\c_{n_j}\} \subseteq\{ \c_{l} \ ~|~ \ l\in N\}$ is linearly independent, which contradicts with Eq. \eqref{eqn:relationship}. Hence, in the first place, any dictionary with a co-basis not attaining the lexicographical minimum order is not the end of pivoting by the Zero rule. {It implies that the collection of terminal dictionaries of any hyperplane arrangement is single}, which is the co-basis with the lexicographical minimum order in the arrangement.

\rightline{$\Box$}

{\textbf{Remark 2:} It can be seen that the proposed Zero rule is independent of the column $x_g$ which corresponds to the constant terms in the dictionary. This means that when the constant terms in the dictionary change, the selection made by the Zero rule remains unchanged. According to our definition of the dictionary, these constant terms originate from the right-hand side constants of the hyperplanes, namely $b_1, \cdots, b_n$. Changes in these constants geometrically represent the translation of the hyperplanes, but in the dictionary representation, they only affect the column $x_g$. Therefore, we can conclude that the Zero pivot is translation invariant.}

\textbf{Connection to an open question}. Propositions \ref{prop: d-steps conv} and \ref{prop:convergence} might be related to an open question: "Is it possible to find a polynomial pivot rule for linear programming, \textit{i.e.}, a pivot rule which bounds the number of pivot steps by a polynomial function of the number of variables and the number of inequalities, or to prove that such a pivot rule does not exist". This question is important \cite{terlaky1993pivot}, motivated by analyzing the complexity of the simplex method in linear programming. The diameter of a polytope provides a lower bound on the number of steps needed by the simplex method, although it is actually weaker than the so-called Hirsch conjecture in polyhedral theory. The Hirsch conjecture, formulated by W. M. Hirsch in 1957 in \cite{klee1987d} and reported in the 1963 book of Dantzig \cite{dantzig1963linear}, states that the edge-vertex graph of an $n$-facet polytope in $d$-dimensional Euclidean space has a diameter no more than $n-d$. This means that any two vertices of a polytope must be connected by a path of length at most $n-d$. {Although there are many studies on the conjecture and diameter of a polytope \cite{Blanchard2020OnTL, michini2014hirsch, Rispoli1998ABO}},
the conjecture is shown to be true for (0,1)-polytopes \cite{naddef1989hirsch} and generally false by a counterexample presented in \cite{Santos2010ACT} which introduces a 43-dimensional polytope of 86 facets with a diameter exceeding 43.  

However, this counterexample does not affect the analysis of the simplex method. It includes the construction of a pivot rule starting from outside the feasible region with an upper bound of $n-d$. The recently developed Facet pivot rule exemplifies the latter possibility. It can achieve an optimal dictionary within at most $n-d$ pivot steps and enters the feasible region from the outside in the final step \cite{yang2021facet,yang2022facet}. Nevertheless, its crucial aspect lies in ensuring dual feasibility at each step, making it difficult to aid the construction of a pivot rule that allows movement between adjacent vertices within the feasible region \cite{yang2021facet,yang2022facet}. In contrast, the Zero rule does not impose stringent conditions in each pivot step. Therefore, our rule, bounded by a linear function of the dimension, may shed some light on the open question. The Zero rule is based on the dictionary itself and goes without a ratio test. Thus, it does not inherently traverse between adjacent vertices or reflect the number of edges between vertices. However, if a ratio test could be incorporated or primal feasibility could be maintained within the Zero rule, it might serve as a polynomial pivoting rule on polytopes.



\section{VE Algorithm Using the Zero Rule}
\label{sec4}
In the above, we introduce the Zero rule. {Here, we describe the algorithm proposed by Avis and Fukuda \cite{avis1991pivoting}, introduce our algorithm using the Zero rule, and conduct a detailed complexity analysis.}


\subsection{Description of the \texttt{AF} algorithm \cite{avis1991pivoting}} 

{The \texttt{AF} algorithm \cite{avis1991pivoting} is based on the following principle: within the entire hyperplane arrangement, starting from any given dictionary, the Criss-Cross pivot uniquely provides a finite-length path to an optimal dictionary. All such paths collectively form a graph rooted in this optimal dictionary. Therefore, we can begin from the optimal dictionary and apply valid reverse Criss-Cross pivots in lexicographic order, which corresponds to a depth-first search over the graph. Each dictionary is counted once. By searching every optimal dictionary, we can obtain all dictionaries within the hyperplane arrangement. Meanwhile, to address the situation where a vertex corresponds to multiple dictionaries, we only need to conduct a test when the algorithm discovers a new dictionary to ensure that among all dictionaries associated with the same vertex, only one meets the test criteria, which will then be output or recorded.}

\begin{figure}[htbp]
\vspace{-0.3cm}
    \centering
    \includegraphics[width=0.6\linewidth]{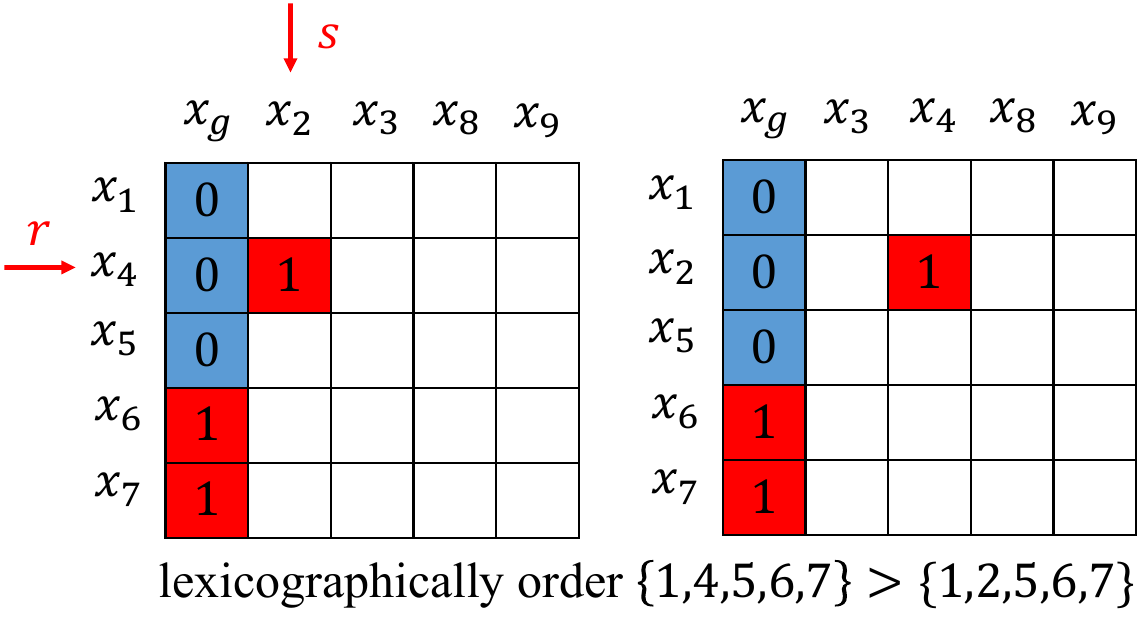}
    \caption{An example of the lex-min test. LHS: fails the test. RHS: passes the test.}
    \label{fig:lex-min test}
    \vspace{-0.2cm}
\end{figure}

T{hus, the algorithm can be divided into three steps: first, construct an initial optimal dictionary from the given vertex; perform a reverse search related to the Criss-Cross rule for each optimal dictionary, outputting all dictionaries that meet the testing criteria; after the reverse search, generate another optimal dictionaries based on the currect one and then go second step. If there is no more optimal dictionaries, the algorithm terminate.} The key aspects of this algorithm are the testing criteria and the method used to conduct the reverse search. As Figure \ref{fig:lex-min test} shows, the testing criterion established by Avis and Fukuda involves selecting the dictionary whose basis is lexicographically minimum among all those corresponding to the same basic solution. {For a non-degenerate dictionary associated with a non-degenerate vertex, this dictionary is the only one associated with that vertex; therefore, no need to compare. Regarding the degenerate dictionary associated with a degenerate vertex, multiple degenerate dictionaries corresponding to the same vertex may exist; thus, we need to apply the following standard for selection.}

\begin{prop}\label{prop: lex-min test}
    Let $B$ be a basis for a degenerate dictionary $\x_B=\bar{A}\x_N$. $B$ is not lexicographically minimum for the corresponding basic solution if and only if there exists $r\in B_{\neq f}$ and $s\in N_{\neq g}$ such that $r>s, \bar{a}_{rg}=0$ and $\bar{a}_{rs}\neq0$.~\\
\end{prop}

\begin{figure}[htbp]
    \centering
    \includegraphics[width=0.95\linewidth]{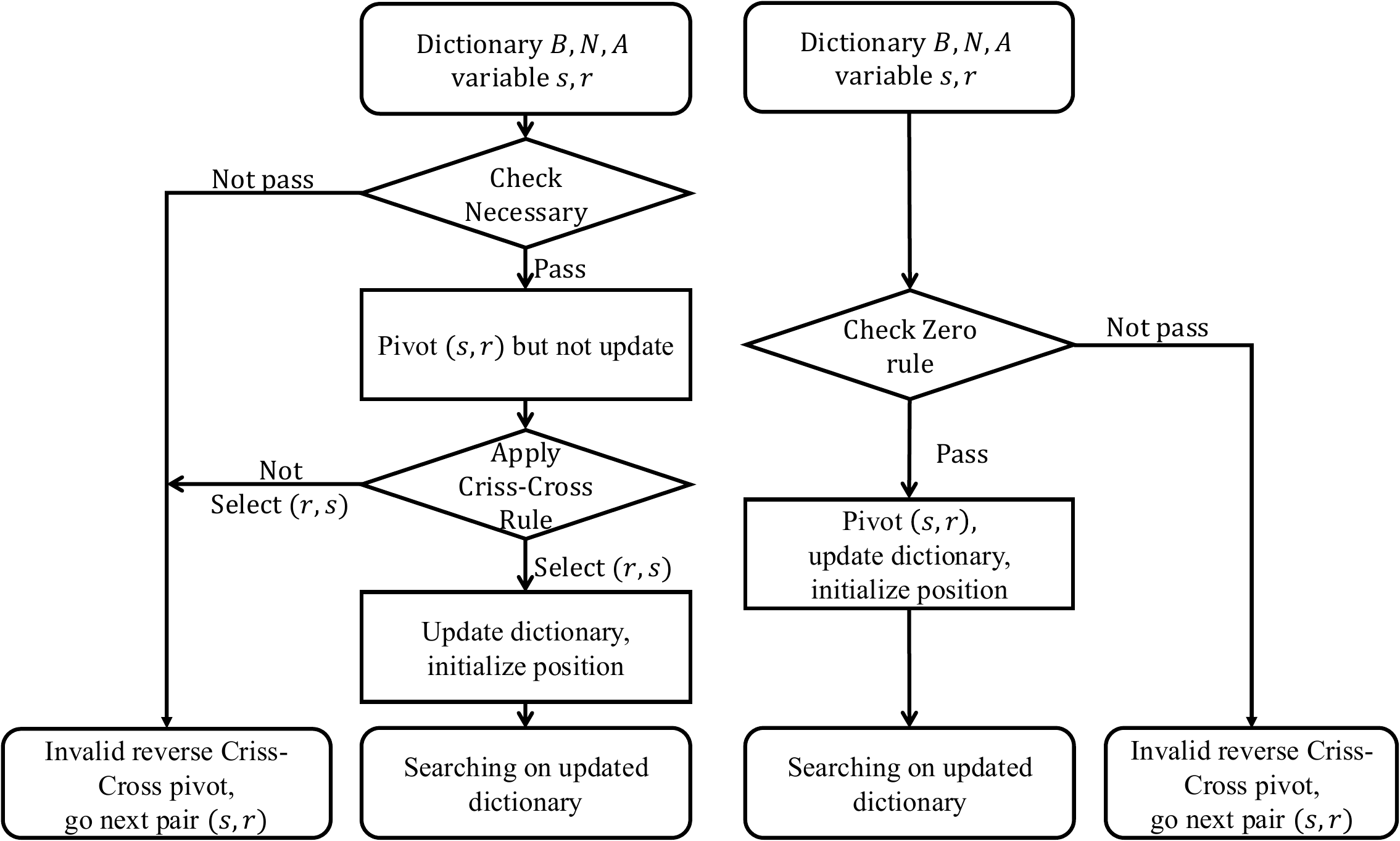}
    \vspace{0.2cm}
    \caption{The process of checking and performing valid reverse pivots. LHS: Avis and Fukuda \cite{avis1991pivoting}; RHS: Our proposed algorithm.}
    \label{fig:reverse_check}
\end{figure}

{As for the depth-first search mentioned earlier, since our goal is to find the paths generated by the Criss-Cross pivot, it is necessary to determine whether a pivot $(s, r)$, where $x_s \in B$ and $x_r \in N$, on a dictionary is a valid reverse Criss-Cross pivot. To do this, the pivot is performed, but we do not use it to update the dictionary that the currently investigated, and the Criss-Cross rule is applied to the resulting dictionary. If the selection yields $(r, s)$, then the pivot $(s, r)$ on the currently investigated dictionary is a valid reverse Criss-Cross pivot. We update the currently investigated dictionary to the computed dictionary and initialize the checking position. Otherwise, do not update and move directly to the next position to repeat the checking process.} Additionally, Avis and Fukuda \cite{avis1991pivoting} provide a necessary condition to test before performing the pivot, which helps reduce some unnecessary computations during the process. The LHS of Figure \ref{fig:reverse_check} illustrates the whole process.

If every position in the current dictionary has been checked, they then apply the Criss-Cross rule for selection and perform a Criss-Cross pivot, followed by updating the dictionary and the search position. If the Criss-Cross rule does not select a proper position, the reverse search concerning one optimal dictionary is finished, thereby proceeding with another optimal dictionary.

{\textbf{Remark 3}. It can be observed that determining whether a pivot is a valid reverse Criss-Cross pivot is quite complex, and sometimes results in $\mathcal{O}(d(n-d)+n)=\mathcal{O}(nd)$ redundant operations. To address this issue, we can construct an if-and-only-if condition to enhance the \texttt{AF} algorithm \cite{avis1991pivoting}. Because this is beyond the scope of the proposed Zero rule and our algorithm, we put this if-and-only-if condition and its detailed analysis and comparison in the \hyperref[analysis_iff]{Appendix}, which may be of independent interest to some readers. We denote the algorithm using this if-and-only-if condition as the \texttt{Enhanced AF}. Although its theoretical operational complexity and space complexity are the same as the original \texttt{AF}, namely $\mathcal{O}(n^2dv)$ and $\mathcal{O}(nd)$, in practice, the redundant operations in their algorithm introduce several lower-order terms, which can only be neglected given a large $d$ and a larger $n$. Furthermore, even if they can be ignored in theoretical analysis, they still have a considerable impact during actual execution, whereas the \texttt{Enhanced AF} {does not have} these terms. In all, we think that the \texttt{Enhanced AF} is more efficient.}

\subsection{Vertex Enumeration with the Zero Rule and Advantage Analysis}

{Now, we integrate the Zero rule into the algorithmic framework of Avis and Fukuda \cite{avis1991pivoting} for the vertex enumeration. First, we note that the core of the \texttt{AF} algorithm is to have a rule to generate a path that leads to an optimal dictionary of the entire hyperplane arrangement. This means that any other pivoting rule including the Zero rule that generates a path towards an optimal dictionary of the entire hyperplane arrangement could also seamlessly fit this framework. Algorithm \ref{alg:VE_alg_Zero_rule} shows the proposed new algorithm with the Zero rule.}


\textbf{Advantage analysis and detailed formulation}. Since our focus is on the VE of an arrangement, we only compare with the Criss-Cross rule and Jensen’s general relaxed recursion that are applicable to this problem. The advantages of the Zero rule are as follows:

Jensen's general relaxed recursion requires partitioning a dictionary and analyzing it case by case. Its termination dictionary could be primal feasible, dual feasible, or optimal \cite{jensen1985coloring, terlaky1993pivot}, which makes it tedious for VE. As for the Criss-Cross rule, although it also needs to discuss two cases, it has only one type of termination dictionary, namely the optimal dictionary. This makes it better than Jensen's general relaxed recursion. The Criss-Cross rule is the exact one used in the \texttt{AF} algorithm \cite{avis1991pivoting}.

However, compared to the Zero rule, the Criss-Cross rule is not ideal in both initialization and the reverse search. During the initialization process, the Zero rule is more convenient in generating all dictionaries where the pivoting rule cannot select a position. By Proposition \ref{prop:convergence}, through our way of initialization, the resultant dictionary with co-basis $\{g,1,2,\cdots,d\}$ is always the unique terminal dictionary. Therefore, after applying the Zero rule to the algorithm, the overall enumeration process can be simplified. More favorably, since the terminal dictionary is the dictionary with the smallest lexicographic order co-basis among all dictionaries in the hyperplane arrangement, it can be obtained simply by renaming the hyperplanes. Let the first dictionary we construct satisfy $N = \{g, 1, 2, \ldots, d\}$ $B = \{d+1, \ldots, n\}$; this automatically be the unique terminal dictionary. While using the Criss-Cross rules requires setting an objective function and generating all optimal dictionaries, which adds an extra step to the initialization process.

{When it comes to the reverse search, the Zero rule is also more favorable. According to Proposition \ref{prop: d-steps conv}, the number of pivot steps in the Zero rule is bounded by $d$, whereas that in the Criss-Cross rule is exponential \cite{Roos1990AnEE, terlaky1993pivot}. We introduce the concept of layers into the algorithm. We assign $\mathtt{layer}=0$ to the terminal dictionary, and when a valid reverse Zero pivot is performed, we update $\mathtt{layer} = \mathtt{layer} + 1$; conversely, when a Zero pivot is performed, we update $\mathtt{layer} = \mathtt{layer} - 1$. Since the number of pivot steps in the Zero rule is bounded by $d$, there are no dictionaries at the $(d+1)$-th layer, meaning that there is no need to do the valid reverse Zero pivots at the $d$-th layer, which can save the cost of checking the valid reverse.}

{Fortunately, we can adopt the idea from Avis and Fukuda \cite{avis1991pivoting} of "computing but not updating”. For dictionaries at the $(d-1)$-th layer, if the pivot $(s, r)$ is a valid reverse Zero pivot, we compute the pivot $(s, r)$ but do not update the dictionary being examined. After determining whether to output the resultant dictionary, we move directly to the next position, thereby bypassing all redundant computations on dictionaries at the $d$-th layer. However, the Criss-Cross rule requires examining every dictionary, which incurs more computational costs.}

\begin{figure}[htbp]
    \centering
    \includegraphics[width=0.6\linewidth]{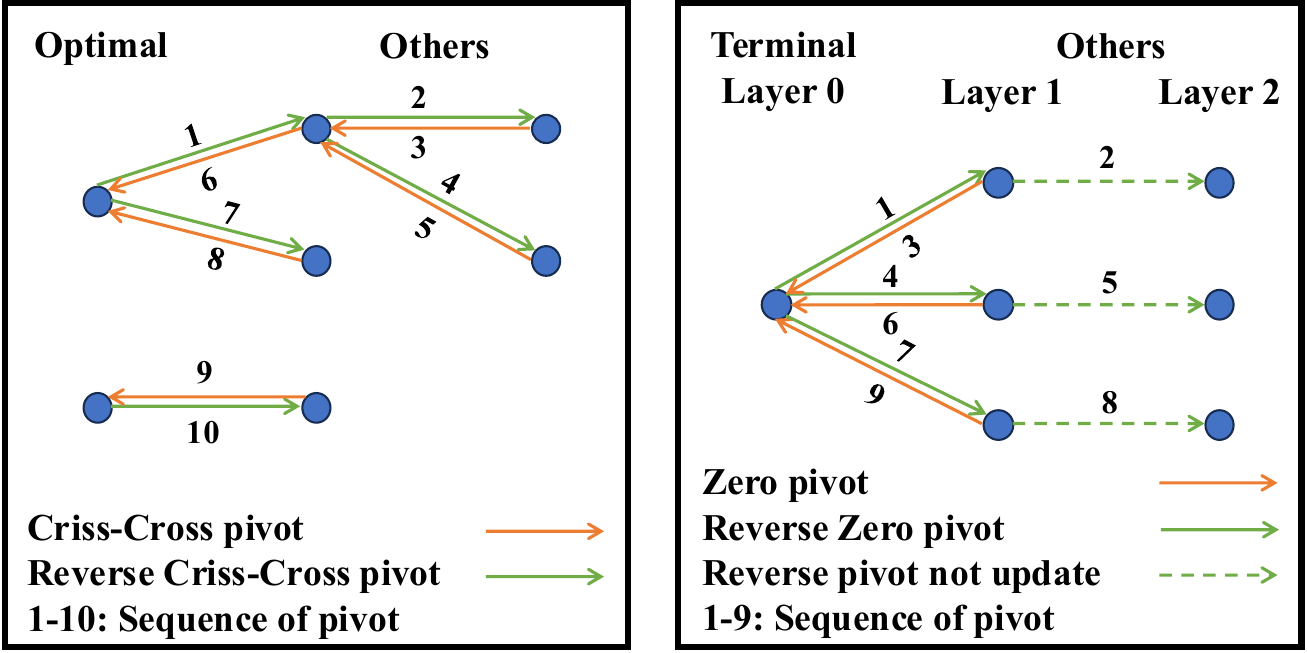}
    \caption{A graphic explanation of the \texttt{Enhanced AF} (LHS) and ours (RHS).}
    \label{fig:search_explanation}
    \vspace{-0.5cm}
\end{figure}

A graphic explanation of our algorithm is on the RHS of Figure \ref{fig:search_explanation}, and the pseudocode for it is presented in Algorithm \ref{alg:VE_alg_Zero_rule}.

\begin{algorithm}[htbp]
    \caption{Ours, doing the reverse search using the Zero rule}
    \begin{algorithmic} 
        \Function{$\mathtt{Search}(B,N,\bar{A})$}{}\\
            \hspace{0.5cm}$i=1,j=2,\mathtt{layer}=0$;\hfill{\% Use $i, j$ represent the $i$-th row, $j$-th column in $\bar{A}$.}
            \If {$\mathtt{lex-min}(B,N,\bar{A})$==1}\\
                \hspace{1cm}Print $B$;
            \EndIf
            \While{$j\leq \mathtt{length}(N)+1$}
                \If{$j\leq \mathtt{length}(N)$}
                    \If{$\mathtt{reverse}(B,N,\bar{A},i,j)==1$ and $\mathtt{layer}<d-1$} \\
                        \hspace{2.1cm}$[B,N,\bar{A}]=\mathtt{pivot}(B,N,\bar{A},i,j)$;\hfill{ \% Compute and update.} 
                        \If {$\mathtt{lex-min}(B,N,\bar{A})==1$}\\
                            \hspace{2.6cm}Print $B$;
                        \EndIf\\
                        \hspace{2.1cm}$i=1,j=2,\mathtt{layer}=\mathtt{layer}+1;$\hfill{\% Initialize the position, go next layer.}
                    \Else \textbf{if} ~ {$\mathtt{reverse}(B,N,\bar{A},i,j)==1$ and $\mathtt{layer}==d-1$} \\
                    \hspace{2.1cm}$[tempB,tempN,temp\bar{A}]=\mathtt{pivot}(B,N,\bar{A},i,j)$;\hfill{ \% Not update.} 
                        \If {$\mathtt{lex-min}(tempB,tempN,temp\bar{A})==1$}\\
                            \hspace{2.6cm}Print $tempB$;
                        \EndIf\\
                        \hspace{2.1cm}$[i,j]=\mathtt{increment}(i,j)$;\hfill{\% Go next position.}\\
                    \hspace{1.6cm}\textbf{else}\\
                        \hspace{2.1cm}$[i,j]=\mathtt{increment}(i,j)$;\hfill{\% Invalid, go next position.}
                    \EndIf
                \Else\hfill{\% Each position on this dictionary has been checked.}\\
                    \hspace{1.6cm}$[i,j]=\mathtt{select}(B,N,\bar{A})$;\hfill{\% Use Zero rule to select a position.}
                    \If{both $i$ and $j$ nonempty}\hfill{\% The Zero rule find a position.}\\
                        \hspace{2.1cm}$[B,N,\bar{A}]=\mathtt{pivot}(B,N,\bar{A},i,j)$;\hfill{\% Back to parent dictionary.}\\
                        \hspace{2.1cm}$[i,j]=\mathtt{find}(i,j)$;\hfill{\% Find the corresponding position.}\\
                        \hspace{2.1cm}$[i,j]=\mathtt{increment}(i,j)$;\\
                        \hspace{2.1cm}$\mathtt{layer}=\mathtt{layer}-1$;
                    \Else\hfill{\% The Zero rule does not find a position.}\\
                        \hspace{2.1cm}Break;\hfill{\% The function $Search$ terminated.}
                    \EndIf
                \EndIf
            \EndWhile
        \EndFunction \vspace{0.2cm}
        
        \noindent\textbf{function} $\mathtt{lex-min}(B,N,\bar{A})$\hfill{\% Test if the current basis output.}\\
        \noindent\textbf{function} $\mathtt{reverse
        }(B,N,\bar{A},i,j)$\hfill{\% Test if pivot$(B(i),N(j))$ a valid reverse pivot.}\\
        \noindent\textbf{function} $[B,N,\bar{A}]=\mathtt{pivot}(B,N,\bar{A},i,j)$\hfill{\% Compute the new dictionary.}\\
        \noindent\textbf{function} $[i,j]=\mathtt{increment}(i,j)$\hfill{\% Go next position on the dictionary.} \\
        \noindent\textbf{function} $[i,j]=\mathtt{select}(i,j)$\hfill{\% Use the Zero rule to find a position.}\\
        \noindent\textbf{function} $[i,j]=\mathtt{find}(i,j)$\hfill{\% Find corresponding position after Zero pivot.}
    \end{algorithmic}
    \label{alg:VE_alg_Zero_rule}
\end{algorithm}

\subsection{Complexity Analysis} 
\label{sec:ca}
Given a hyperplane arrangement in $\mathbb{R}^d$ consisting of $n$ distinct hyperplanes with $v$ vertices, there are $\mathcal{O}(v)$ dictionaries. We note that applying the proposed Zero pivot takes at most $d$ iterations. Therefore, we denote the number of dictionaries that require exactly $d$ iterations as $\mathcal{O}(v_d)$. Note that when calculating a new dictionary, one needs $(n-d)(d+1)=\mathcal{O}(nd)$ calculations. Testing the lex-min also needs $\mathcal{O}(nd)$ calculations.

\subsubsection{Complexity of The \texttt{AF} and The \texttt{Enhanced AF} Algorithm} 
{The computational complexity of the \texttt{Enhanced AF} Algorithm consist of the lex-min test on each dictionary, checking and performing the valid reverse pivot, and performing the Criss-Cross pivot. The total operation required for checking and performing the valid reverse pivot can be found in \hyperref[analysis_iff]{Appendix}, Eq. \eqref{eq:alg_reformulate_RHS}, which is $\mathcal{O}(n^2dv)$. The lex-min tests and the Criss-Cross pivot together require at most $\mathcal{O}(ndv+(n + nd)v) = \mathcal{O}(ndv)$ operations. Therefore, its computational complexity is $$\mathcal{O}(ndv + n^2dv) = \mathcal{O}(n^2dv).$$
As for the space complexity, since there is no additional storage, we only need to store those variables generated during the computation, which is in the order of $\mathcal{O}(nd)$.


{Before proceeding with the analysis of the \texttt{AF} algorithm, it is important to note that the algorithm aims to maintain a balanced efficiency among all hyperplane arrangements. Therefore, we need to analyze the general case. Thus, we can assume that the probability of any position fulfilling the examination condition in each dictionary is $0.5$, regardless of the necessary condition in the original \texttt{AF} algorithm or the if-and-only-if condition used in the \texttt{Enhanced AF} algorithm. To facilitate the analysis in the following paragraph, we state the necessary condition in \cite{avis1991pivoting} as follows:}

\begin{prop}\label{prop:cc_neces}
    If $(s, r), s\in B_{\neq f}, r\in N_{g}$, is a valid reverse criss-cross pivot for a dictionary $x_B = Ax_N$, then either:
    \begin{align*}
        (a)&~\bar{a}_{sg}>0,\bar{a}_{sr}>0,\bar{a}_{sj}\geq0 \ for \ j\in N_{\neq g},j<s,\\
        (b)&~\bar{a}_{fr}<0,\bar{a}_{sr}<0,\bar{a}_{ir}\leq0 \ for  \ i\in B_{\neq f},i<r.
    \end{align*}
\end{prop}

Note that the \texttt{AF} algorithm requires processing each position given every dictionary exactly once. This means that checking and performing the valid reverse pivot in this algorithm requires $\mathcal{O}(ndv)$ times processing in the LHS of Figure \ref{fig:reverse_check}, where there are $\mathcal{O}(v)$ valid reverse pivots.

{We can categorize $\mathcal{O}(ndv)$ pivot check into three types: valid reverse pivots, those meet the necessary condition but are invalid, and those that fail the necessary condition. For the first two types, their cost is $\mathcal{O}(nd + n + n) = \mathcal{O}(nd)$. Observing Proposition \ref{prop:cc_iff} in \hyperref[analysis_iff]{Appendix} and the necessary condition, they require checking at most $2(2d + 2(n - d)) = 4n$ positions and $d + (n - d) = n$ positions, respectively. Thus, the second type occurs with a probability of $2^{-n + d - 1} + 2^{-d - 1} - 2^{-2n - 1} - 2^{-2n - 1} = 2^{-n + d - 1} + 2^{-d - 1} - 2^{-2n}$. Therefore, the total cost for the first two types combined is 
\begin{align*}
    &\mathcal{O}\Big(ndv + (2^{-n + d - 1} + 2^{-d - 1} - 2^{-2n})n^2d^2v\Big)\\
    =&\mathcal{O}\Big(ndv + (2^{-n+d}+2^{-d })n^2d^2v\Big).
\end{align*}
The third type incurs a cost of 
\begin{align*}
    &\mathcal{O}\Big((ndv - v - (2^{-n + d - 1} + 2^{-d - 1} - 2^{-2n})ndv)n\Big)\\
    =&\mathcal{O}\Big((nd - 1 - (2^{-n+d} + 2^{-d})nd)nv\Big)\\
    =&\mathcal{O}\Big((1 - 2^{-n+d} - 2^{-d})n^2dv\Big)
\end{align*}
operations. {Thus, the total operations in the \texttt{AF} algorithm requires for checking and performing valid reverse pivots are}
\begin{equation}
    \begin{aligned}
        &\mathcal{O}\Big(ndv + (2^{-n+d}+2^{-d })n^2d^2v+(1 - 2^{-n+d} - 2^{-d})n^2dv\Big)\\
        =&\mathcal{O}\Big(ndv + (2^{-n+d}+2^{-d })(d-1)n^2dv+n^2dv\Big)\\
        =&\mathcal{O}\Big((2^{-n+d}+2^{-d })n^2d^2v+n^2dv\Big).
    \end{aligned}\label{eq:alg_avis_LHS}
\end{equation}}

{In addition to the computations in Eq. \eqref{eq:alg_avis_LHS}, the consuming of those lex-min tests and those Criss-Cross pivots also require at most $\mathcal{O}(ndv+(n + nd)v) = \mathcal{O}(ndv)$ operations. Therefore, its total computational complexity should be}
\begin{align*}
\mathcal{O}\Big((2^{-n+d}+2^{-d})n^2d^2v+n^2dv+ndv\Big)=\mathcal{O}(n^2dv).
\end{align*}

{Similar with the analysis in the \texttt{Enhanced AF} algorithm, it needs to apply the reverse search to each optimal dictionary, and thus the space complexity should be $\mathcal{O}(nd)$.}


{Although this computational complexity is the same as that of \texttt{Enhanced AF}, we have omitted many lower-order terms compared to \texttt{Enhanced AF}. Moreover, the $(2^{-n+d}+2^{-d})n^2d^2v$ could be negligible only when the dimension is sufficiently large and the number of hyperplanes should be even larger, which is quite difficult to achieve in practice. Therefore, in our empirical experiments, the \texttt{Enhanced AF} performs better.}

\subsubsection{Complexity of The \texttt{Moss} Algorithm} 

Now, let us analyze the computational complexity of the method from \texttt{Moss} \cite{moss2012basis}. In each dictionary, one needs to compute the ratio of the column of $x_g$ vs every column and then select the extreme value. This step incurs a cost of $\mathcal{O}(d(n-d)v)= \mathcal{O}(ndv)$. Subsequently, during the comparison with discovered dictionaries, their co-basis needs to be checked. If there are $k$ recorded dictionaries, then the comparison requires $\mathcal{O}(kd)$ for each candidate. When generating the adjacent co-basis about a fixed dictionary, one may find at most $nd$ candidates from it. Thus, comparing each candidate takes $\mathcal{O}(knd^2)$ computations. Since there are a total of $\mathcal{O}(v)$ dictionaries, the total amount of computations is $\mathcal{O}(\sum_{k=1}^v(knd^2))=\mathcal{O}(nd^2v^2)$. It is important to note that one only needs to compute those dictionaries whose co-basis are not duplicated, which costs $\mathcal{O}(ndv)$ in computing dictionaries and determining if they pass the lex-min test. Consequently, the complexity of the method from \texttt{Moss} \cite{moss2012basis} in VE is $\mathcal{O}(ndv+ nd^2v^2+ndv)= \mathcal{O}(nd^2v^2)$. 

As for the space complexity, since the \texttt{Moss} algorithm needs to store and study each dictionary in the arrangement, the space complexity should be $\mathcal{O}(ndv)$.

\subsubsection{Complexity of Our Algorithm} 
It is important to note that each dictionary needs to be computed at most twice, namely through the valid reverse Zero pivot to enter it and through the Zero pivot to leave it. Therefore, the total cost of computing those dictionaries should be $\mathcal{O}(ndv)$. For each entry of a dictionary, testing a valid reverse Zero pivot requires at most $\frac{1}{2}(n-d+n)d$ examinations, which is of $\mathcal{O}(nd)$. While $\mathcal{O}(v-v_d)$ dictionaries need to be tested for validity, $\mathcal{O}(v_d)$ dictionaries do not. Thus, the total complexity incurred by testing validity is $\mathcal{O}(n^2d^2(v-v_d))$. In addition, each selection in the Zero rule also brings $\mathcal{O}(nd)$ comparisons, and we need to apply the selection on $\mathcal{O}(v-v_d)$ dictionaries. Together with the operations conducted in the lex-min test, the complexity of our algorithm is approximately
\begin{equation}
    \begin{aligned}
    \mathcal{O}\Big(ndv+n^2d^2(v-v_d)+ndv+ndv\Big)=\mathcal{O}\Big(n^2d^2(v-v_d)+ndv_d\Big)
    \end{aligned}
    \label{eqn:complexity}
    \end{equation}
{When it comes to the space complexity, since there is no additional storage, we only need to store those variables generated during the computation, which can take the space of at most two dictionaries, in the order of $\mathcal{O}(nd)$. We organize the complexity of our algorithm into a theorem:}

\begin{theorem}
    {Given a hyperplane arrangement in $\mathbb{R}^d$ consisting of $n$ distinct hyperplanes and $v$ vertices, let $\mathcal{O}(v_d)$ denote the number of dictionaries that require exactly $d$ steps of the Zero pivot. Then, the computational and space complexity of our algorithm is $\mathcal{O}(n^2d^2(v-v_d) + ndv_d)$ and $\mathcal{O}(nd)$, respectively.}
\end{theorem}

The ratio of $v_d$ vs $v$ has a significant impact on the actual complexity. Below we use two typical classes of hyperplane arrangements to illustrate the effect of $v_d$: one is a special case that consists of a lot of parallel hyperplanes; the other is a simple arrangement. Before estimation, we first characterize what kinds of dictionaries lie in the $d$-th layer.

\begin{prop}[Disjoint layer d]
    Assume that the termination dictionary about the Zero rule is with the co-basis $N_0$, then those dictionaries with $N\cap N_0=\emptyset$ are in the $d$-th layer.
\end{prop}

\noindent\textit{Proof: } If those dictionaries with $N\cap N_0=\emptyset$ are not in the $d$-th layer, after repeating applying the Zero rule on the dictionary to termination, there remains at least one index in $N$ such that it is also in the co-basis of termination, which is $N_0$. This causes a contradiction.

\rightline{$\Box$}

\begin{figure}[htbp]
    
    \centering
    \includegraphics[width=0.6\textwidth]{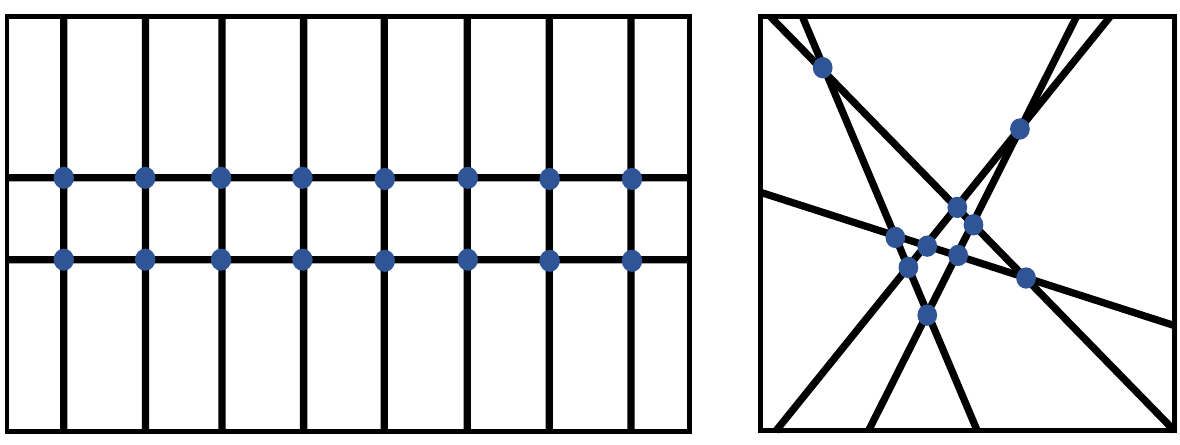}
    \caption{LHS: a sequence of (hyper)cubes. RHS: a simple arrangement.}
    \label{fig:two type arrangements} 
    \vspace{-0.7cm}
\end{figure}

In an arrangement with a lot of parallel hyperplanes, where hyperplanes form boundaries of hypercubes arranged along a specific coordinate axis, such as $y_1, y_2, \ldots, y_d = 0, 1, y_1 = 2, y_1 = 3, y_1 = 4, \ldots,y_1=n-2d+1$. As shown in the LHS of Figure \ref{fig:two type arrangements}, $v_d =2^{-d}v$, then Eq. \eqref{eqn:complexity} can be simplified to
$$\mathcal{O}(n^2d^2(v-v_d)+ndv_d)=\mathcal{O}((1-2^{-d})n^2d^2v+2^{-d}ndv)=\mathcal{O}(n^2d^2v).$$
If one stacks multiple layers of hypercubes on top of each other along this line, the proportion of $v_d$ will increase. Consequently, as the number of hypercubes increases, the complexity for this class of hyperplane arrangements will significantly decrease.

In the case of a simple arrangement (see the RHS of Figure \ref{fig:two type arrangements}, $v=\binom{n}{d}, v_d\geq\binom{n-d}{d}$), then  Eq.\eqref{eqn:complexity} can be reduced to
\begin{align*}
    \mathcal{O}(n^2d^2(v-v_d)+ndv_d)&\leq\mathcal{O}\left(n^2d^2\left(\binom{n}{d}-\binom{n-d}{d}\right)+nd\binom{n-d}{d}\right)\\
    &=\mathcal{O}\left(n^2d^2\left(\frac{n^d-(n-d)^d}{d!}+\mathcal{O}(n^{d-2})\right) +nd\binom{n-d}{d}\right)\\
    &=\mathcal{O}\left(n^2d^2\left(\frac{d^2 n^{d-1}}{d!}+\mathcal{O}(n^{d-2})\right) +nd\binom{n-d}{d}\right)\\
    &=\mathcal{O}\left(n^2d^2\left(d\frac{n^{d-1}}{(d-1)!}+\mathcal{O}(n^{d-2})\right) +nd\binom{n-d}{d}\right)\\
    &=\mathcal{O}\left(n^2d^3\binom{n}{d-1}+nd\binom{n-d}{d}\right)\\
    &=\mathcal{O}\left(nd^4\binom{n}{d}+nd\binom{n-d}{d}\right)\\
    &=\mathcal{O}\left(nd^4\binom{n}{d}\right)=\mathcal{O}(nd^4v).
\end{align*}

In fact, the simple arrangement is a quite general case. If a hyperplane arrangement is randomly generated, it is very likely to be a simple arrangement, and therefore the complexity of our algorithm is significantly smaller than its counterparts. We reasonably contend that our algorithm is state-of-the-art in the problem of VE of an arrangement.

The lower bound of the computational complexity of the reversal VE algorithm is $\mathcal{O}(ndv)$, since at least given each vertex, one needs to take $\mathcal{O}(nd)$ to obtain another dictionary. Our algorithm achieves this lower bound at $v_d$ vertices. This is because the valid reverse examination is no longer needed for vertices in the $d$-th layer. The Zero rule can guarantee the termination in $d$ steps. Can we further reduce the computational complexity of the proposed algorithm? It is difficult. To do so, one may seek a new pivot rule that either increases the portion of dictionaries that waive a valid reverse examination or reduces the number of examinations per valid reverse test. In the first situation, the potential superior rule needs to have pivot steps in $\mathcal{O}(\log(d))$ or $\mathcal{O}(1)$ to terminate. However, this requirement is impossible, as there is always a certain dictionary that terminates in at least $d$ steps \textit{i.e.}, those dictionaries discussed in Proposition \ref{sec5}. In the second situation, reducing the examination means less information from a dictionary but more pivot steps. Given these challenges, the proposed Zero rule is a good choice in terms of achieving the lower bound over at least $v_d$ vertices.

\section{Toy Example}\label{sec5}

\begin{wrapfigure}{r}{0.3\textwidth}
    \vspace{-0.4cm}
    \centering
    \includegraphics[width=0.3\textwidth]{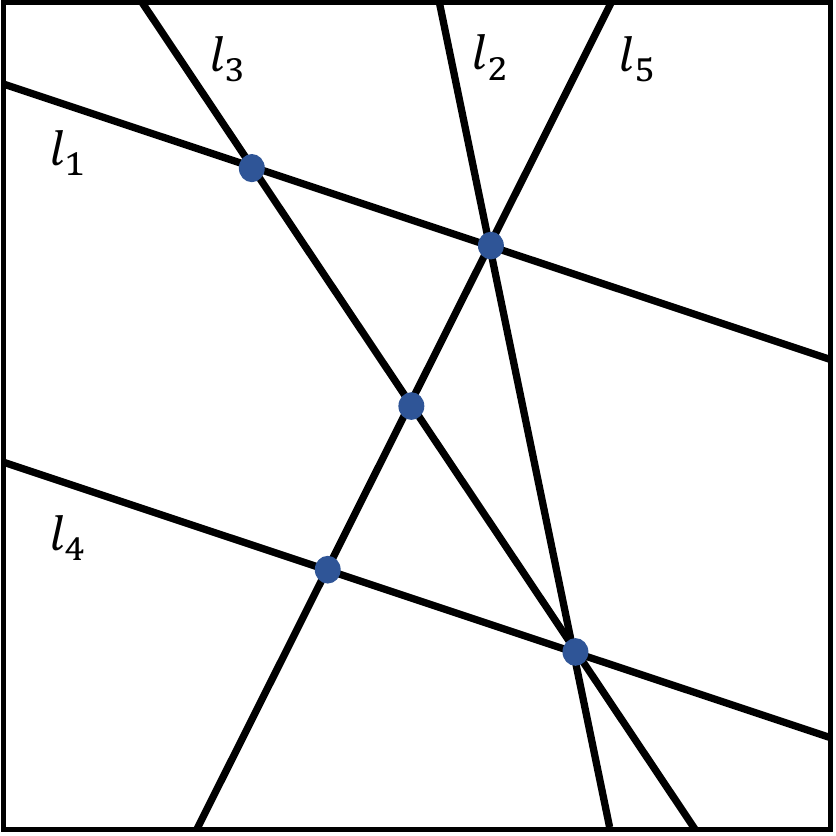}
    \caption{Representative example.}
    \label{fig:The representative example}
    \vspace{-0.8cm}
\end{wrapfigure}

In this section, we present a simple yet illustrative example to demonstrate the advantages of employing the Zero rule. This example highlights the limitations of employing the Criss-Cross rule: firstly, the Criss-Cross pivot does not guarantee convergence within $d$ steps. Secondly, regardless of the chosen objective function, there may exist multiple optimal dictionaries, resulting in multiple search trees. It is noteworthy that the example here is not particularly special, since the weights and biases of the hyperplane are generated randomly. 


Let us consider the arrangement in Figure \ref{fig:The representative example}, which has 5 hyperplanes, and the weights and biases of each hyperplane are 
\begin{equation}
\begin{cases}
    &\c_1=(1,3), ~~~~~ b_1=4 \\
    &\c_2=(5,1), ~~~~~ b_2=5 \\
    &\c_3=(3,2), ~~~~~ b_3=2 \\
    &\c_4=(-1,-3), ~b_4=1 \\
    &\c_5=(-2,1), ~~~b_5=-\frac{1}{2},
\end{cases}
\end{equation}
and slack variables are $x_1, ...,x_5$, which are expressed in the following equations:
{\begin{equation*}
\begin{cases}
    &x_1=4x_g-y_1-3y_2\\
    &x_2=5x_g-5y_1-y_2\\
    &x_3=2x_g-3y_1-2y_2\\
    &x_4=x_g+y_1+3y_2\\
    &x_5=-\frac{1}{2}x_g+2y_1-y_2
\end{cases}
\end{equation*}}

\textbf{Ours}: Taking $x_1$ and $x_2$ as co-basis and adding the variable $x_g$, we obtain the following initial dictionary:
\begin{equation}
\begin{cases}
    x_3&=-\frac{5}{2}x_g+\frac{1}{2}x_1+\frac{1}{2}x_2\\
    x_4&=5x_g-x_1+0x_2\\
    x_5&=0x_g +\frac{1}{2}x_1-\frac{1}{2}x_2.
\end{cases}
\label{eqn:layer0 dic}
\end{equation}

The dictionary with $N=\{g,1,2\}$ is the lexicographically minimal among all co-basis in this arrangement. We assign it with $\mathtt{layer}=0$. {As for the lex-min test, for convenience, denote this dictionary as $(B^0, N^0, \bar{A}^0)$. Since there exists a pair of variables $(5, 1)$, where $5 \in B^0$ and $1 \in N^0$, such that $\bar{a}_{3g}^0 = 0$ and $\bar{a}_{32}^0 = \frac{1}{2} \neq 0$, $B^0$ is not lexicographically minimal for the corresponding basic solution. Therefore, it did not pass the lex-min test, and we do not need to print its basis.}

Next, we examine each entry sequentially, pivot $(3,1)$ on the dictionary \eqref{eqn:layer0 dic} is not a valid reverse Zero pivot, while $(4,1)$ is. Note that the layer of dictionary \eqref{eqn:layer0 dic} is $0$, which is less than $d-1=1$, then we pivot $(4,1)$ and obtain a dictionary in $\mathtt{layer}=0+1=1$, and update the dictionary as below: 
\begin{equation}
\begin{cases}
    x_1&=5x_g+ 0x_2-x_4\\
    x_3&=0x_g +\frac{1}{2}x_2-\frac{1}{2}x_4\\
    x_5&=\frac{5}{2}x_g-\frac{1}{2}x_2-\frac{1}{2}x_4.
\end{cases}
\label{eqn:layer1 dic}
\end{equation}

{Now we do the lex-min test for this dictionary. Denote this dictionary as $(B^1, N^1, \bar{A}^1)$, we see there exists a pair of variables $(3, 2)$, where $3 \in B^1$ and $2 \in N^1$, such that $\bar{a}_{3g}^1 = 0$ and $\bar{a}_{32}^1 = \frac{1}{2} \neq 0$. According to Proposition \ref{prop: lex-min test}, $B^1$ is not lexicographically minimum for the corresponding basic solution. }Therefore, it does not pass the lex-min test, and we do not print its basis. 

After the lex-min test, we resume the reverse search in this dictionary. The pivot $(1,2)$ is invalid, and the pivot $(3,2)$ is valid. Note that the layer of dictionary \eqref{eqn:layer1 dic} is $1$, which is equal to $d-1=1$, then we pivot $(3,2)$ and obtain the following dictionary in $\mathtt{layer}=1+1=2$ temporarily but do not update the dictionary being examined:
\begin{equation}
\begin{cases}
    x_1&=5x_g+0x_3-x_4\\
    x_2&=0x_g +2x_3+x_4\\
    x_5&=\frac{5}{2}x_g-x_3 - x_4.
\end{cases}
\end{equation}
{Denote the dictionary by $(B^{temp},N^{temp},\bar{A}^{temp})$, the only variable $r\in B^{temp}$ such that $\bar{a}^{temp}_{rg} = 0$ is $2$. Since $\forall s\in N_{\neq g}, 2<s$, Proposition \ref{prop: lex-min test} implies this $B^{temp}$ attain the lexicographically minimum for the corresponding basic solution; therefore, we output its basis $\{1,2,5\}$.}

Consequently, 
we continue searching on the dictionary being examined, which is the dictionary \eqref{eqn:layer1 dic}. After checking all entries on this dictionary, all remaining entries are invalid reverses, the dictionary \eqref{eqn:layer1 dic} is finished. Then, we can apply the Zero rule on it to go to the dictionary \eqref{eqn:layer0 dic} and so forth. The LHS of Figure \ref{fig:example result spanning} presents a spanning tree that encompasses all dictionaries and their orders in our algorithm.

\begin{figure}[htbp]
    \centering
    \includegraphics[width=0.75\linewidth]{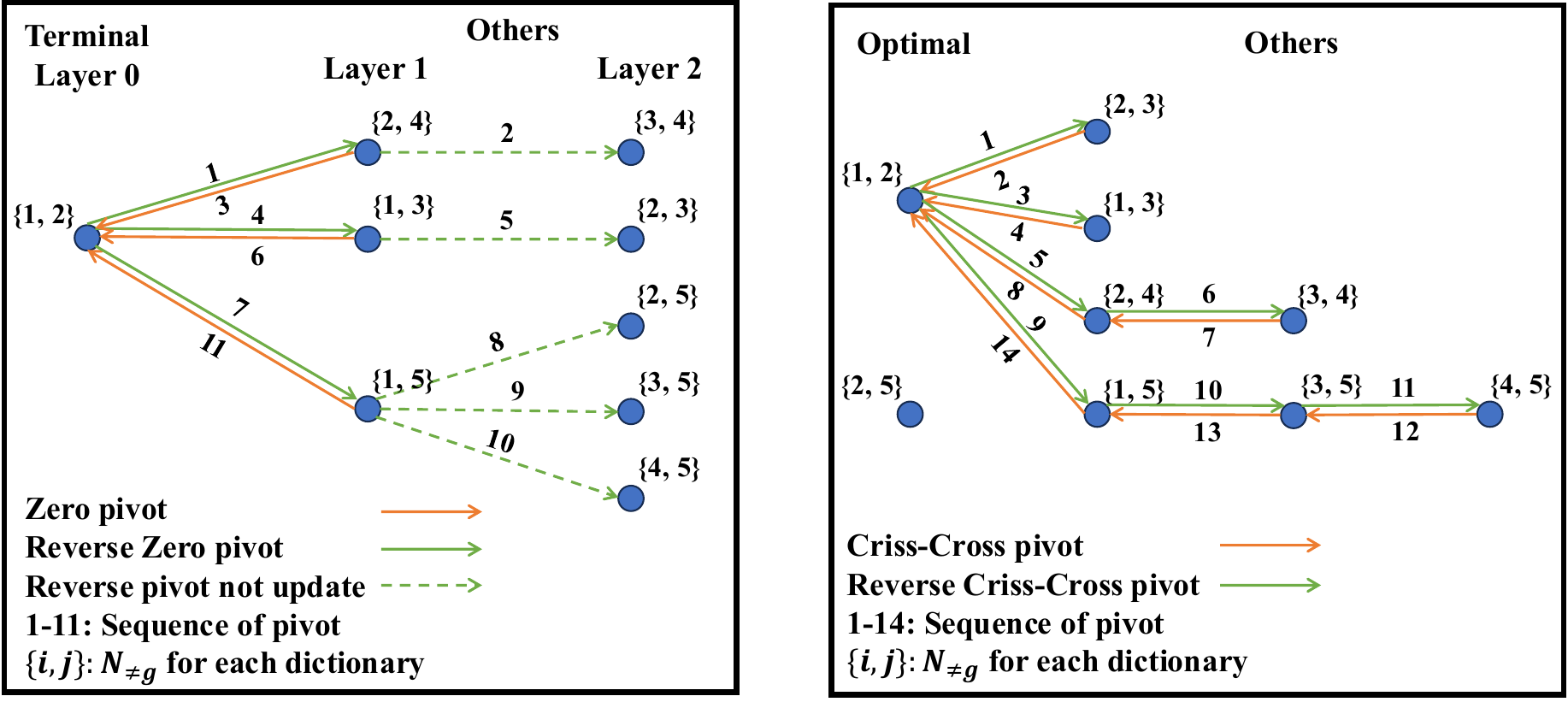}
    \caption{Spanning trees that connect all vertices obtained by the VE algorithm. LHS: Ours. RHS: \texttt{Enhanced AF}.}
    \label{fig:example result spanning}
    \vspace{-0.5cm}
\end{figure}

\textbf{Enhanced AF}: Now, we utilize the \texttt{Enhanced AF} algorithm to solve this small example. The first step is to construct an optimal dictionary with a proper objective function. According to \cite{avis1991pivoting}, we set the initial dictionary as the dictionary \eqref{eqn:ini dic cc}, which is adding a row of objective function to dictionary \eqref{eqn:layer0 dic}. 
\begin{equation}
\begin{cases}
    x_3&=-\frac{5}{2}x_g+\frac{1}{2}x_1+\frac{1}{2}x_2\\
    x_4&=5x_g-x_1+0x_2\\
    x_5&=0x_g +\frac{1}{2}x_1-\frac{1}{2}x_2\\
    x_f&=0x_g-x_1-x_2.
\end{cases}
\label{eqn:ini dic cc}
\end{equation}
Meanwhile, another optimal dictionary in the hyperplane arrangement is
\begin{equation}
\begin{cases}
    x_1&=0x_g+x_2+2x_5\\
    x_3&=-\frac{5}{2}x_g+x_2+x_5\\
    x_4&=5x_g -x_2-2x_5\\
    x_f&=0x_g-x_2-2x_5.
\end{cases}
\label{eqn:ini dic cc other}
\end{equation}

The algorithm needs to to apply the reverse search about the Criss-Cross rule on them, respectively. Let us start from the dictionary \eqref{eqn:ini dic cc}.

Similar with the analysis to the dictionary \eqref{eqn:layer0 dic}, this dictionary does not pass the lex-min test, and the pivot $(3,1)$ is a valid reverse Criss-Cross pivot. Pivoting $(3,1)$ can obtain the following dictionary:
\begin{equation}
\begin{cases}
    x_1&=5x_g-x_2-2x_3\\
    x_4&=0x_g-x_2+2x_3\\
    x_5&=\frac{5}{2}x_g -x_2-x_3\\
    x_f&=-5x_g+0x_2+2x_3.
\end{cases}
\label{eqn:dic cc {2,3}}
\end{equation}

{Denote this dictionary by $(B^2,N^2,\bar{A}^2)$, then there exists a pair of variables $(4,2), 4\in B_{\neq f}^2, 2\in N^2$ such that $\bar{a}^2_{4g}=0$ and $\bar{a}^2_{42}=-1\neq 0$. It implies that the dictionary fails the lex-min test.}

After carefully checking all entries of this dictionary, all remaining entries are invalid reverses, and the study on dictionary \eqref{eqn:dic cc {2,3}} is finished. We can apply the Criss-Cross rule on it and go to the dictionary \eqref{eqn:ini dic cc}. Next, we continue the test for the pair of variables after $(3,1)$ and so forth.

On the other hand, we apply the reverse search with the Criss-Cross rule on the dictionary \eqref{eqn:ini dic cc other}. There is no valid reverse Criss-Cross pivot on it; therefore, the dictionary \eqref{eqn:ini dic cc other} is an isolated dictionary. Thus, the entire VE algorithm terminated. The RHS of Figure \ref{fig:example result spanning} presents a spanning tree that encompasses all dictionaries and their orders by \texttt{Enhanced AF}. 

\section{Systematic Experiments}
\label{sec6}

Here, we demonstrate the superiority of our algorithm through systematic experiments. We compare our algorithm with \texttt{Enhanced AF}, \texttt{AF}, and \texttt{Moss} on the following four types of hyperplane arrangements: an arrangement containing $2d$ hyperplanes that enclose a unit cube in $\mathbb{R}^d$; $2d+1$ hyperplanes that enclose a truncated unit cube in $\mathbb{R}^d$, \textit{i.e.}, adding the hyperplane $y_1+y_2+\cdots+y_d=1.5$ with the unit cube; 125 different hyperplane arrangements randomly generated by MATLAB; and a family of arrangements that consist of all possible arrangements generated by a three-layer ReLU neural network with weights randomly generated. For each example, we list the number of vertices found by four algorithms and the run time. We denote the time spent exceeding 10,000 seconds as $>$9,999e purpose of the first two arrangements is that the ground truth is known, which can assist us to examine if the compared algorithms can enumerate all vertices upon completion of the algorithm. The third is to show the comparison in the general case. The last one is to show the utility of the proposed algorithm in artificial networks. All our examples and code are publicly available\footnote{https://github.com/Github-DongZelin/Examples-and-Command-of-the-algorithm} for readers' free download and use.

\textbf{Unit hypercubes}: Tables \ref{table:cubes_number} and \ref{table:cubes_time} presents different algorithms' performance over unit hypercubes ranging from $\mathbb{R}^2$ to $\mathbb{R}^{8}$. 

\begin{table}[htbp]
    \renewcommand{\arraystretch}{1.5}
    \centering
    \begin{tabular}{|p{2.2cm}<{\centering}|p{1.2cm}<{\centering}|p{1.2cm}<{\centering}|p{1.2cm}<{\centering}|p{1.2cm}<{\centering}|p{1.2cm}<{\centering}|p{1.2cm}<{\centering}|p{1.2cm}<{\centering}|p{1.2cm}<{\centering}|p{1.2cm}<{\centering}|p{1.2cm}<{\centering}|}
         \hline
          Method & $\mathbb{R}^2$ & $\mathbb{R}^{3}$ & $\mathbb{R}^4$ & $\mathbb{R}^5$ & $\mathbb{R}^6$ & $\mathbb{R}^7$ & $\mathbb{R}^8$  \\
         \hline
         Truth & 4 & 8 &  16 & 32 & 64 & 128 & 256 \\\hline
         {Ours} & \multirow{4}{*}{4} & \multirow{4}{*}{8} & \multirow{4}{*}{16} & \multirow{4}{*}{32} & \multirow{4}{*}{64} & \multirow{4}{*}{128} & \multirow{4}{*}{256} \\\hhline{|-|~|~|~|~|~|~|~|~|}
         \texttt{AF} \cite{avis1991pivoting} &  &  &  &  &  &  & \\ \hhline{|-|~|~|~|~|~|~|~|~|}
         \texttt{Enhanced AF} &  &  &  &  &  &  & \\ \hhline{|-|~|~|~|~|~|~|~|~|}
         \texttt{Moss} \cite{moss2012basis} &  &  &  &  &  &  & \\ \hline
    \end{tabular}
    \caption{The number of vertices found by four algorithms over unit hypercubes from $\mathbb{R}^2$ to $\mathbb{R}^{8}$.}
    \label{table:cubes_number}
    \vspace{-0.2cm}
\end{table}

\begin{table}[htbp]
    \renewcommand{\arraystretch}{1.5}
    \centering
    \begin{tabular}{|p{2.2cm}<{\centering}|p{1.2cm}<{\centering}|p{1.2cm}<{\centering}|p{1.2cm}<{\centering}|p{1.2cm}<{\centering}|p{1.2cm}<{\centering}|p{1.2cm}<{\centering}|p{1.2cm}<{\centering}|p{1.2cm}<{\centering}|p{1.2cm}<{\centering}|p{1.2cm}<{\centering}|}
         \hline
          Method & $\mathbb{R}^2$ & $\mathbb{R}^{3}$ & $\mathbb{R}^4$ & $\mathbb{R}^5$ & $\mathbb{R}^6$ & $\mathbb{R}^7$ & $\mathbb{R}^8$  \\
         \hline
         Ours
          & 0.2201 & 0.6535 & 1.9113 & 4.8336 & 12.264 & 30.774 & 73.478 \\\hline
         \texttt{AF} \cite{avis1991pivoting}
          & 0.5261 & 1.5073 & 4.0578 & 10.476 & 26.228 & 64.643 & 150.81 \\\hline
        \texttt{Enhanced AF}
          & 0.4023 & 1.0381 & 2.7246 & 6.9270 & 17.287 & 43.149 & 105.63 \\\hline
         \texttt{Moss} \cite{moss2012basis}
          & 0.3520 & 0.8673 & 2.5098 & 8.0306 & 30.945 & 129.95 & 628.63 \\\hline
    \end{tabular}
    \caption{The time consuming by four algorithms over unit hypercubes from $\mathbb{R}^2$ to $\mathbb{R}^{8}$.}
    \label{table:cubes_time}
    \vspace{-0.2cm}
\end{table}

\textbf{Unit hypercubes with a cone being cut}: Tables \ref{table:cubecones_time} and \ref{table:cubecones_number} presents different algorithms' performance over this arrangement ranging from $\mathbb{R}^2$ to $\mathbb{R}^{8}$.

\begin{table}[htbp]
    \renewcommand{\arraystretch}{1.5}
    \centering
    \begin{tabular}{|p{2.2cm}<{\centering}|p{1.2cm}<{\centering}|p{1.2cm}<{\centering}|p{1.2cm}<{\centering}|p{1.2cm}<{\centering}|p{1.2cm}<{\centering}|p{1.2cm}<{\centering}|p{1.2cm}<{\centering}|p{1.2cm}<{\centering}|p{1.2cm}<{\centering}|p{1.2cm}<{\centering}|}
         \hline
         Method & $\mathbb{R}^2$ & $\mathbb{R}^{3}$ & $\mathbb{R}^4$ & $\mathbb{R}^5$ & $\mathbb{R}^6$ & $\mathbb{R}^7$ & $\mathbb{R}^8$
         \\\hline
         Truth & 8 & 20 & 48 & 112 & 256 & 576 & 1280\\\hline
         {Ours} & \multirow{4}{*}{8} & \multirow{4}{*}{20} & \multirow{4}{*}{48} & \multirow{4}{*}{112} & \multirow{4}{*}{256} & \multirow{4}{*}{576} & \multirow{4}{*}{1280} \\\hhline{|-|~|~|~|~|~|~|~|~|}
         \texttt{AF} \cite{avis1991pivoting} &  &  &  &  &  &  & \\ \hhline{|-|~|~|~|~|~|~|~|~|}
         \texttt{Enhanced AF} &  &  &  &  &  &  & \\ \hhline{|-|~|~|~|~|~|~|~|~|}
         \texttt{Moss} \cite{moss2012basis} &  &  &  &  &  &  & \\ \hline
    \end{tabular}
    \caption{The number of vertices found by four algorithms over unit hypercubes with a cone being cut from $\mathbb{R}^2$ to $\mathbb{R}^{8}$.}
    \label{table:cubecones_number}
    \vspace{-0.2cm}
\end{table}

\begin{table}[htbp]
    \renewcommand{\arraystretch}{1.5}
    \centering
    \begin{tabular}{|p{2.2cm}<{\centering}|p{1.2cm}<{\centering}|p{1.2cm}<{\centering}|p{1.2cm}<{\centering}|p{1.2cm}<{\centering}|p{1.2cm}<{\centering}|p{1.2cm}<{\centering}|p{1.2cm}<{\centering}|p{1.2cm}<{\centering}|p{1.2cm}<{\centering}|p{1.2cm}<{\centering}|}
         \hline
         Method & $\mathbb{R}^2$ & $\mathbb{R}^{3}$ & $\mathbb{R}^4$ & $\mathbb{R}^5$ & $\mathbb{R}^6$ & $\mathbb{R}^7$ & $\mathbb{R}^8$  \\\hline
         {Ours}
         & 0.4817 & 1.5661 & 5.1182 & 15.792 & 46.351 & 129.48 & 378.05 \\\hline
         \texttt{AF} \cite{avis1991pivoting}
         & 1.5312 & 4.9375 & 15.173 & 45.272 & 124.24 & 340.59 & 912.61 \\\hline
         \texttt{Enhanced AF}
         & 1.0575 & 2.8983 & 8.7185 & 25.462 & 71.393 & 195.76 & 547.49 \\\hline
         \texttt{Moss} \cite{moss2012basis}
         & 0.9378 & 3.2863 & 16.060 & 96.004 & 647.23 & 4314.6 & $>$9999 \\\hline
    \end{tabular}
    \caption{The time consuming by four algorithms over unit hypercubes with a cone being cut from $\mathbb{R}^2$ to $\mathbb{R}^{8}$.}
    \label{table:cubecones_time}
    \vspace{-0.2cm}
\end{table}

\textbf{Random arrangement}. By taking random integers as coefficients of arrangements, a total of 25 cases are generated by varying from $\mathbb{R}^2$ to $\mathbb{R}^6$ and from 8 to 16 hyperplanes. Each case has five runs to ensure the reliability of the results. Table \ref{table:ave-ver} displays the average number of vertices found by three algorithms. Table \ref{table:ave-ti} shows the average computation time.

\begin{table}[htbp]
    \vspace{-1cm}
    \renewcommand{\arraystretch}{1.5}
    \centering
    \begin{tabular}{|p{0.8cm}<{\centering}|p{2.2cm}<{\centering}|p{1.2cm}<{\centering}|p{1.2cm}<{\centering}|p{1.2cm}<{\centering}|p{1.2cm}<{\centering}|p{1.2cm}<{\centering}|} \hline
        $d\backslash n$ & Method  & 8 & 10 & 12 & 14 & 16\\ \hline
        \multirow{4}{*}{$\mathbb{R}^2$} & {Ours} & \multirow{4}{*}{27.2} & \multirow{4}{*}{42} & \multirow{4}{*}{63.8} & \multirow{4}{*}{88.4} & \multirow{4}{*}{110.6} \\\hhline{|~|-|~|~|~|~|~|}
         & \texttt{AF} \cite{avis1991pivoting} &  &  &  &  &  \\ \hhline{|~|-|~|~|~|~|~|}
         & \texttt{Enhanced AF} &  &  &  &  &  \\ \hhline{|~|-|~|~|~|~|~|}
         & \texttt{Moss} \cite{moss2012basis} &  &  &  &  &  \\ \hline
        \multirow{4}{*}{$\mathbb{R}^3$} & {Ours} & \multirow{4}{*}{55.6} & \multirow{4}{*}{119.2} & \multirow{4}{*}{218} & \multirow{4}{*}{362.8} & \multirow{4}{*}{557.6} \\\hhline{|~|-|~|~|~|~|~|}
         & \texttt{AF} \cite{avis1991pivoting} &  &  &  &  &  \\ \hhline{|~|-|~|~|~|~|~|}
         & \texttt{Enhanced AF} &  &  &  &  &  \\ \hhline{|~|-|~|~|~|~|~|}
         & \texttt{Moss} \cite{moss2012basis} &  &  &  &  &  \\ \hline
        \multirow{4}{*}{$\mathbb{R}^4$} & {Ours} & \multirow{4}{*}{69.6} & \multirow{4}{*}{210} & \multirow{4}{*}{494.6} & \multirow{4}{*}{997} & \multirow{4}{*}{1818.8} \\\hhline{|~|-|~|~|~|~|~|}
         & \texttt{AF} \cite{avis1991pivoting} &  &  &  &  &  \\ \hhline{|~|-|~|~|~|~|~|}
         & \texttt{Enhanced AF} &  &  &  &  &  \\ \hhline{|~|-|~|~|~|~|~|}
         & \texttt{Moss} \cite{moss2012basis} &  &  &  &  &  \\ \hline
        \multirow{4}{*}{$\mathbb{R}^5$} & {Ours} & \multirow{4}{*}{56} & \multirow{4}{*}{252} & \multirow{4}{*}{792} & \multirow{4}{*}{2001} & \multirow{4}{*}{4367.6} \\\hhline{|~|-|~|~|~|~|~|}
         & \texttt{AF} \cite{avis1991pivoting} &  &  &  &  &  \\ \hhline{|~|-|~|~|~|~|~|}
         & \texttt{Enhanced AF} &  &  &  &  &  \\ \hhline{|~|-|~|~|~|~|~|}
         & \texttt{Moss} \cite{moss2012basis} &  &  &  &  &  \\ \hline
        \multirow{4}{*}{$\mathbb{R}^6$} & {Ours} & \multirow{4}{*}{28} & \multirow{4}{*}{210} & \multirow{4}{*}{924} & \multirow{4}{*}{3003} & \multirow{4}{*}{8008} \\\hhline{|~|-|~|~|~|~|~|}
         & \texttt{AF} \cite{avis1991pivoting} &  &  &  &  &  \\ \hhline{|~|-|~|~|~|~|~|}
         & \texttt{Enhanced AF} &  &  &  &  &  \\ \hhline{|~|-|~|~|~|~|~|}
         & \texttt{Moss} \cite{moss2012basis} &  &  &  &  &  \\ \hline
    \end{tabular}
    \caption{The average number of vertices found by different algorithms.}
    \label{table:ave-ver}
    \vspace{-0.5cm}
\end{table}

\begin{table}[htbp]
    \renewcommand{\arraystretch}{1.5}
    \centering
    \begin{tabular}{|p{0.8cm}<{\centering}|p{2.2cm}<{\centering}|p{1.2cm}<{\centering}|p{1.2cm}<{\centering}|p{1.2cm}<{\centering}|p{1.2cm}<{\centering}|p{1.2cm}<{\centering}|} \hline
        $d\backslash n$ & Method  & 8 & 10 & 12 & 14 & 16\\ \hline
        \multirow{4}{*}{$\mathbb{R}^2$} & 
        Ours & 1.5421 & 2.5362 & 3.8307  & 5.4066  & 7.0743 \\ \hhline{|~|-|-|-|-|-|-|}
        & \texttt{AF} \cite{avis1991pivoting} & 10.177 & 20.936 & 38.954 & 61.440 & 87.521 \\ \hhline{|~|-|-|-|-|-|-|}
        & \texttt{Enhanced AF} & 4.8326 & 8.9556 & 15.564 & 24.147 & 33.960 \\ \hhline{|~|-|-|-|-|-|-|}
        & \texttt{Moss} \cite{moss2012basis} & 3.7805 & 7.6048 & 14.133 & 24.475 & 37.145 \\ \hline
        \multirow{4}{*}{$\mathbb{R}^3$} & 
        Ours & 4.268 & 9.7160 & 18.370 & 29.983 & 47.327 \\ \hhline{|~|-|-|-|-|-|-|}
        & \texttt{AF} \cite{avis1991pivoting} & 25.643 & 72.388 & 171.80 & 321.27 & 580.23 \\ \hhline{|~|-|-|-|-|-|-|}
        & \texttt{Enhanced AF} & 11.032 & 28.791 & 62.738 & 116.36 & 202.30 \\ \hhline{|~|-|-|-|-|-|-|}
        & \texttt{Moss} \cite{moss2012basis} & 16.572 & 65.205 & 217.43 & 593.14 & 1524.9 \\ \hline
        \multirow{4}{*}{$\mathbb{R}^4$} & 
        Ours & 6.1188 & 20.059 & 50.243 & 112.76 & 369.08 \\ \hhline{|~|-|-|-|-|-|-|}
        & \texttt{AF} \cite{avis1991pivoting} & 32.303 & 135.14 & 416.99 & 1098.1 & 2533.7 \\ \hhline{|~|-|-|-|-|-|-|}
        & \texttt{Enhanced AF} & 13.399 & 51.728 & 147.911 & 378.34 & 955.68 \\ \hhline{|~|-|-|-|-|-|-|}
        & \texttt{Moss} \cite{moss2012basis} & 34.254 & 317.89 & 1960.9 & 8526.4 & $>$9999 \\ \hline
        \multirow{4}{*}{$\mathbb{R}^5$} & 
        Ours & 5.9546 & 30.227 & 103.30 & 401.13 & 1115.9 \\ \hhline{|~|-|-|-|-|-|-|}
        & AF\cite{avis1991pivoting} & 28.430 & 181.00 & 737.39 & 2364.4 & 6653.0 \\ \hhline{|~|-|-|-|-|-|-|}
        & \texttt{Enhanced AF} & 12.506 & 67.570 & 264.56 & 843.64 & 2491.2 \\ \hhline{|~|-|-|-|-|-|-|}
        & \texttt{Moss} \cite{moss2012basis} & 32.466 & 1387.3 & 9736.0 & $>$9999 & $>$9999 \\ \hline
        \multirow{4}{*}{$\mathbb{R}^6$} & 
        Ours & 4.1225 & 26.596 & 138.21 & 857.60 & 2858.1 \\ \hhline{|~|-|-|-|-|-|-|}
        & \texttt{AF} \cite{avis1991pivoting} & 16.046 & 149.52 & 958.72 & 4321.8 & $>$9999 \\ \hhline{|~|-|-|-|-|-|-|}
        & \texttt{Enhanced AF} & 7.1205 & 53.508 & 327.96 & 1340.8 & 4929.3 \\ \hhline{|~|-|-|-|-|-|-|}
        & \texttt{Moss} \cite{moss2012basis} & 13.847 & 732.62 & $>$9999 & $>$9999 & $>$9999 \\ \hline
    \end{tabular}
    \caption{The average computation time consumed by different algorithms.}
    \label{table:ave-ti}
    \vspace{-0.5cm}
\end{table}

\begin{figure}[htbp]
    \centering
    \begin{minipage}{0.45\linewidth}
    \includegraphics[width=\linewidth]{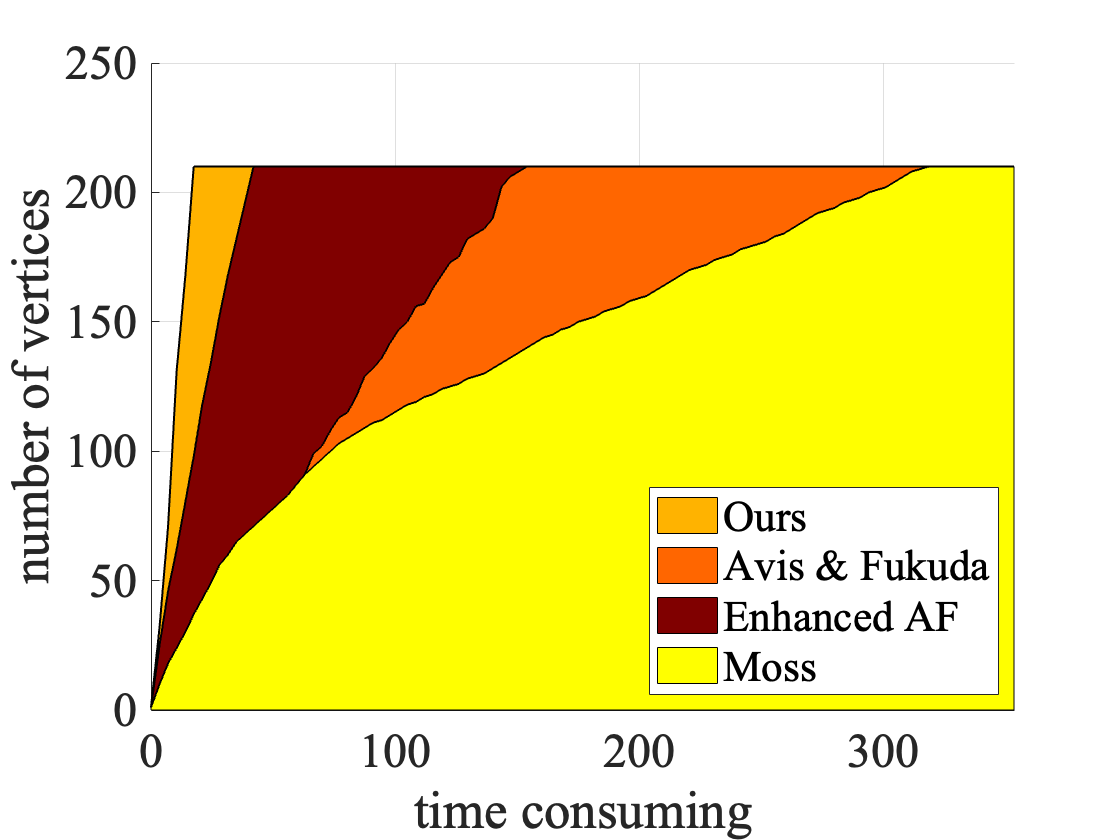}
    \end{minipage}%
    \hspace{0.5cm}
    \begin{minipage}{0.45\linewidth}
    \includegraphics[width=\linewidth]{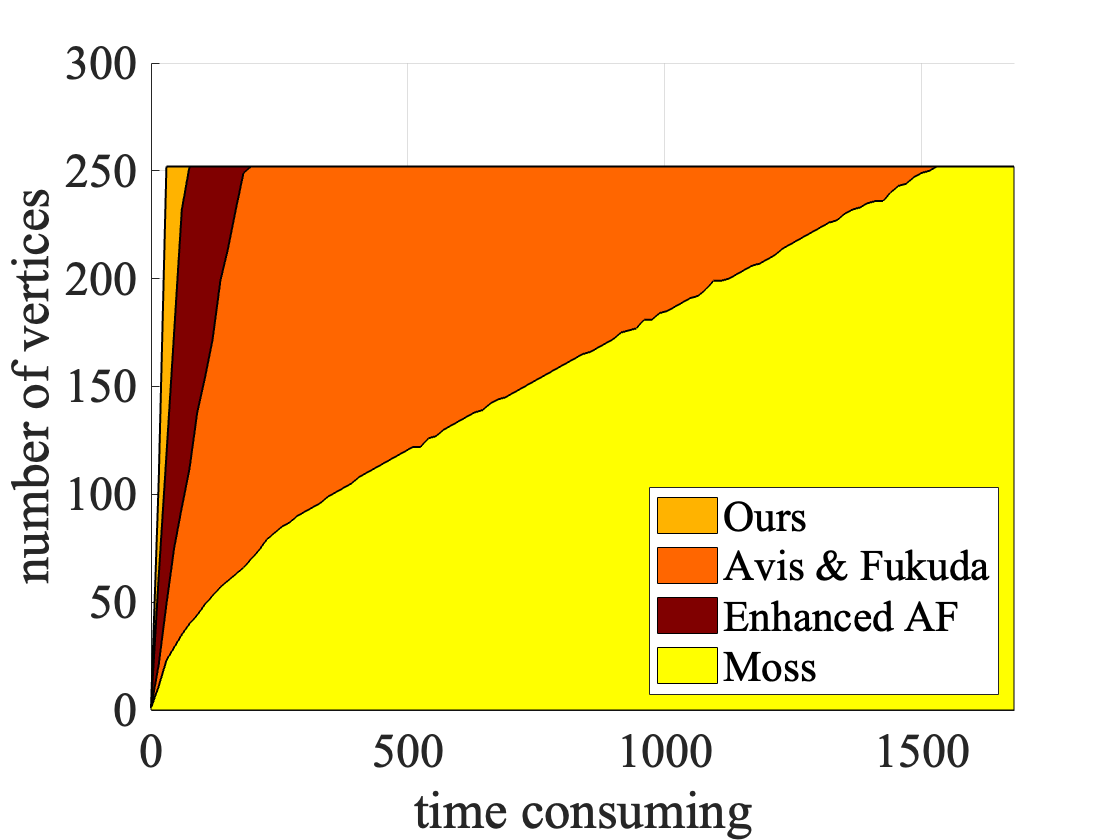}
    \end{minipage}%
    \caption{The number of vertices for an arrangement in $\mathbb{R}^4$(LHS), $\mathbb{R}^5$(RHS) with 10 hyperplanes.}
    \label{fig:rand45}
\end{figure}

Both the \texttt{Enhanced AF} algorithm and ours improve \texttt{Moss} method, Avis and Fukuda's method significantly. One highlight is that reversing the Zero rule is much better than reversing the Criss-Cross rule in the \texttt{Enhanced AF}.

Additionally, from each column in Table \ref{table:ave-ti} and Figures \ref{fig:rand45}, it can be observed that as the number of dimensions and vertices increases, the improvement becomes even more pronounced for any number of hyperplanes.~\\

\textbf{A three-layer ReLU neural network}:
We have built a neural network model with an input layer consisting of $7$ nodes to handle data from $\mathbb{R}^3$. Each node in the input layer receives a $3$-dimensional input. The hidden layer comprises $5$ nodes, each utilizes a weight vector with $7$ dimensions along with a bias term. These layers employ the Rectified Linear Unit (ReLU) activation function to introduce non-linearity into the network. In this network, all weight vectors and bias terms are randomly generated using the rand function in MATLAB. Firstly, we use the auxiliary linear programming to find the partition of the seven neurons in the input layer, resulting in a total of $64$ cells in $\mathbb{R}^3$. As the activation state of each cell is different, $5$ hyperplanes corresponding to the nodes in the hidden are also different. Next, we use algorithms to find out the total number of vertices of these $64$ different hyperplane arrangements with $12=7+5$ hyperplanes in $\mathbb{R}^3$.
We randomly order these arrangements. The partial results are provided in Tables \ref{table:relu-feas/all ver} and \ref{table:relu-time} and the overall result can be found in our link\footnote{https://github.com/Github-DongZelin/Examples-and-Command-of-the-algorithm}, $-1$ in the first arrangement means that none of the seven neurons in the input layer is active.

\begin{table}[htbp]
    \renewcommand{\arraystretch}{1.42}
    \centering
    \begin{tabular}{|p{0.51cm}<{\centering}|p{2.2cm}<{\centering}|p{1.1075cm}<{\centering}|p{1.1075cm}<{\centering}|p{1.1075cm}<{\centering}|p{1.1075cm}<{\centering}|p{1.1075cm}<{\centering}|p{1.1075cm}<{\centering}|} \hline
          & Method & 1 & 4 & 7 & 11 & 14 & 17 \\ \hline
        \multirow{4}{*}{0} & {Ours} & \multirow{4}{*}{-1} & \multirow{4}{*}{185} & \multirow{4}{*}{185} & \multirow{4}{*}{220} & \multirow{4}{*}{220} & \multirow{4}{*}{185} \\\hhline{|~|-|~|~|~|~|~|~|}
         & \texttt{AF} \cite{avis1991pivoting} &  &  &  &  &  &  \\ \hhline{|~|-|~|~|~|~|~|~|}
         & \texttt{Enhanced AF} &  &  &  &  &  & \\ \hhline{|~|-|~|~|~|~|~|~|}
         & \texttt{Moss} \cite{moss2012basis} &  &  &  &  &  & \\ \hline
        \multirow{4}{*}{20} & {Ours} & \multirow{4}{*}{220} & \multirow{4}{*}{220} & \multirow{4}{*}{220} & \multirow{4}{*}{220} & \multirow{4}{*}{220} & \multirow{4}{*}{110} \\\hhline{|~|-|~|~|~|~|~|~|}
         & \texttt{AF} \cite{avis1991pivoting} &  &  &  &  &  &  \\ \hhline{|~|-|~|~|~|~|~|~|}
         & \texttt{Enhanced AF} &  &  &  &  &  & \\ \hhline{|~|-|~|~|~|~|~|~|}
         & \texttt{Moss} \cite{moss2012basis} &  &  &  &  &  & \\ \hline
        \multirow{4}{*}{40} & {Ours} & \multirow{4}{*}{220} & \multirow{4}{*}{220} & \multirow{4}{*}{220} & \multirow{4}{*}{220} & \multirow{4}{*}{220} & \multirow{4}{*}{220} \\\hhline{|~|-|~|~|~|~|~|~|}
         & \texttt{AF} \cite{avis1991pivoting} &  &  &  &  &  &  \\ \hhline{|~|-|~|~|~|~|~|~|}
         & \texttt{Enhanced AF} &  &  &  &  &  & \\ \hhline{|~|-|~|~|~|~|~|~|}
         & \texttt{Moss} \cite{moss2012basis} &  &  &  &  &  & \\ \hline
    \end{tabular}
    \caption{The number of feasible and all vertices in corresponding arrangement found by different algorithms.}
    \label{table:relu-feas/all ver}
    \vspace{0cm}
\end{table}

\begin{table}[htbp]
    \renewcommand{\arraystretch}{1.5}
    \centering
    \begin{tabular}{|p{0.51cm}<{\centering}|p{2.2cm}<{\centering}|p{1.1075cm}<{\centering}|p{1.1075cm}<{\centering}|p{1.1075cm}<{\centering}|p{1.1075cm}<{\centering}|p{1.1075cm}<{\centering}|p{1.1075cm}<{\centering}|} \hline
          & Method & 1 & 4 & 7 & 11 & 14 & 17 \\ \hline
        \multirow{4}{*}{0} 
        & {Ours} & -1 & 16.093 & 16.707 & 18.034 & 18.113 & 16.220 \\\hhline{|~|-|-|-|-|-|-|-|}
        & \texttt{AF} \cite{avis1991pivoting} & -1 & 137.06 & 138.71 & 166.72 & 170.02 & 121.42 \\\hhline{|~|-|-|-|-|-|-|-|}
        & \texttt{Enhanced AF} & -1 & 51.536 & 51.644 & 62.323 & 62.984 & 51.228 \\\hhline{|~|-|-|-|-|-|-|-|}
        & \texttt{Moss} \cite{moss2012basis} & -1 & 146.06 & 148.13 & 218.69 & 223.03 & 149.95 \\\hline
        \multirow{4}{*}{20} 
        & {Ours} & 18.054 & 18.158 & 18.151 & 18.052 & 17.956 & 10.935 \\\hhline{|~|-|-|-|-|-|-|-|}
        & \texttt{AF} \cite{avis1991pivoting} & 161.31 & 171.20 & 173.15 & 166.21 & 175.76 & 70.709 \\\hhline{|~|-|-|-|-|-|-|-|}
        & \texttt{Enhanced AF} & 62.794 & 62.485 & 63.56505 & 62.676 & 62.606 & 30.584 \\\hhline{|~|-|-|-|-|-|-|-|}
        & \texttt{Moss} \cite{moss2012basis} & 220.469 & 223.52 & 227.02 & 227.72 & 227.44 & 52.524 \\\hline
        \multirow{4}{*}{40} 
        & {Ours} & 24.996 & 24.178 & 18.661 & 18.641 & 18.723 & 18.651 \\\hhline{|~|-|-|-|-|-|-|-|}
        & \texttt{AF} \cite{avis1991pivoting} & 198.09 & 237.25 & 170.78 & 179.18 & 179.57 & 178.85 \\\hhline{|~|-|-|-|-|-|-|-|}
        & \texttt{Enhanced AF} & 65.223 & 84.506 & 64.889 & 64.750 & 65.086 & 65.466 \\\hhline{|~|-|-|-|-|-|-|-|}
        & \texttt{Moss} \cite{moss2012basis} & 252.49 & 266.63 & 234.57 & 241.73 & 238.49 & 240.92 \\\hline
    \end{tabular}
    \caption{The average computation time consumed by different algorithms.}
    \label{table:relu-time}
    \vspace{-0.2cm}
\end{table}

From Table \ref{table:relu-feas/all ver}, it can be seen that the vertices in these hyperplane arrangements are highly degenerate. In a simple arrangement formed by 12 hyperplanes in $\mathbb{R}^3$, there should be $220$ vertices. However, several examples do not reach this value. Therefore, this series of examples can test the effectiveness of our method in highly degenerate hyperplane arrangement. 

The following Figures \ref{fig:1737} and Figures \label{fig:2141} demonstrate the process of searching vertices in degenerated and non-degenerated arrangements, respectively, where we switch the order of color coverage in RHS of Figure \ref{fig:1737} for better visualization. It can be seen that, in all cases, ours and the \texttt{Enhanced AF} are the most and the second most efficient, respectively.

\begin{figure}[htbp]
    \centering
    \begin{minipage}{0.45\linewidth}
    \includegraphics[width=\linewidth]{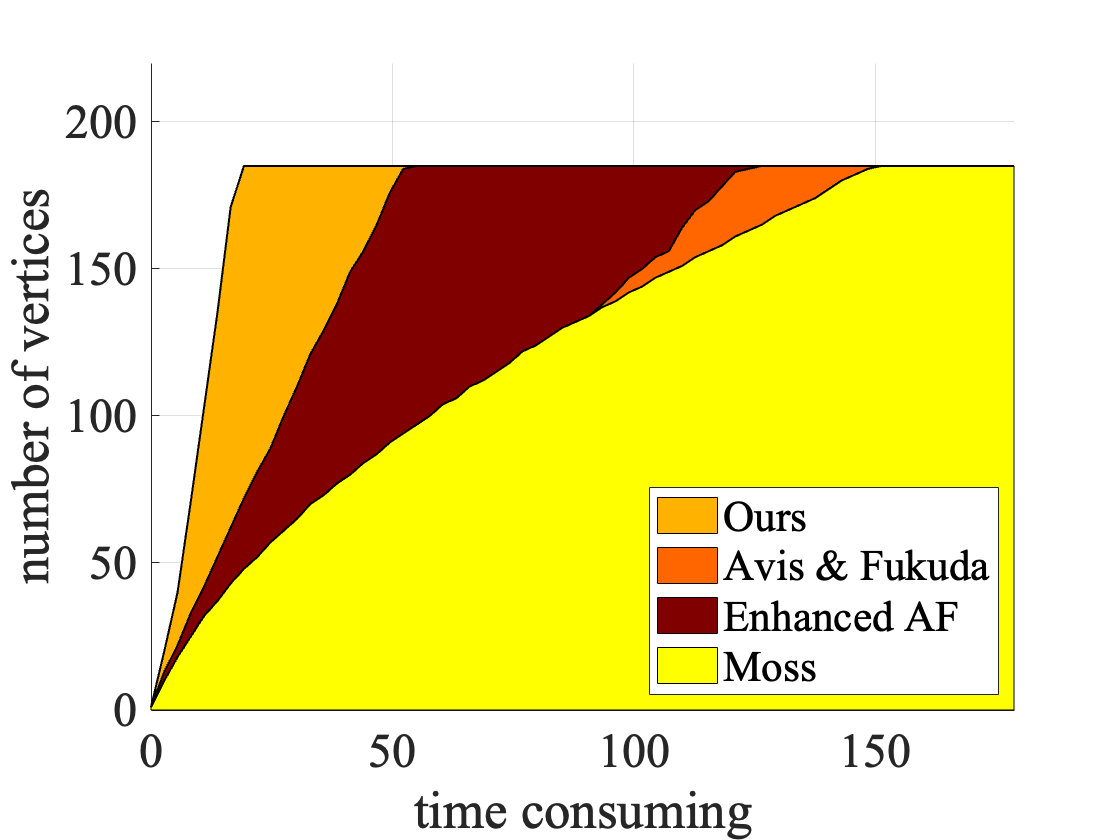}
    \end{minipage}%
    \hspace{0.5cm}
    \begin{minipage}{0.45\linewidth}
    \includegraphics[width=\linewidth]{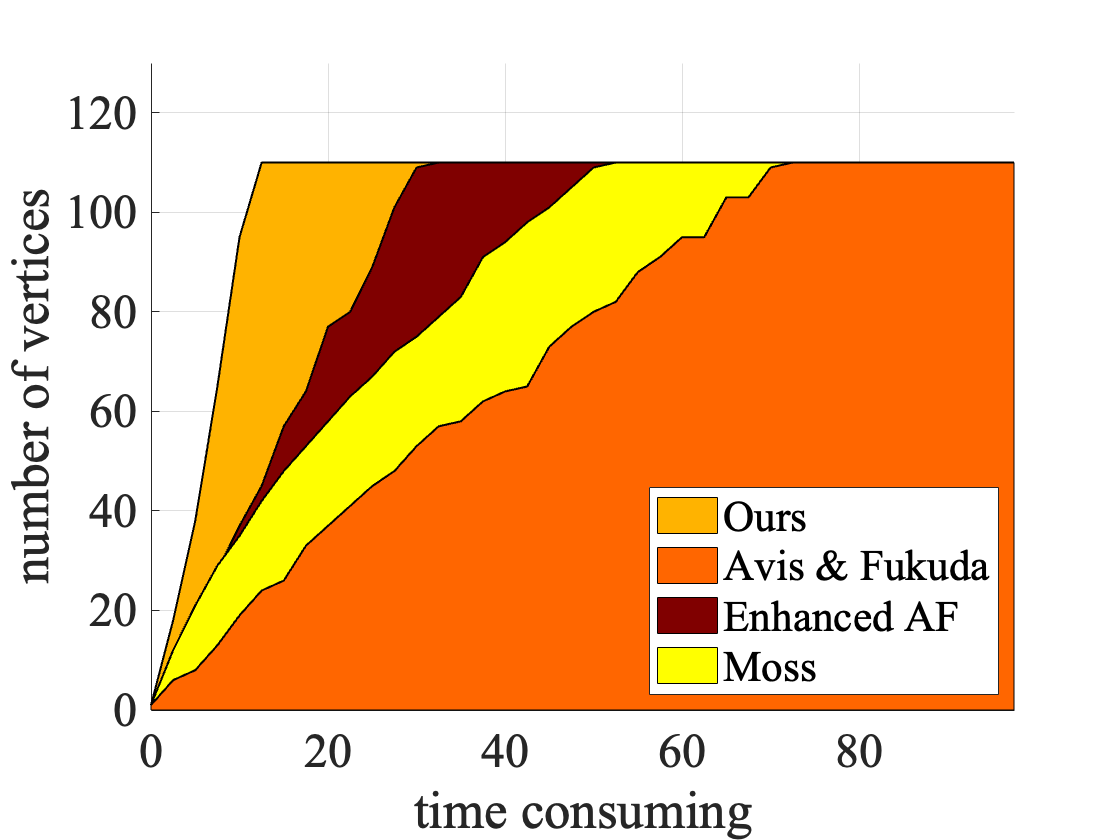}
    \end{minipage}
    \caption{Vertices in 17-th(LHS) and 37-th(RHS) arrangement}
    \label{fig:1737}
    \vspace{-0.4cm}
\end{figure}

\begin{figure}[htbp]
    \centering
    \begin{minipage}{0.45\linewidth}
    \includegraphics[width=\linewidth]{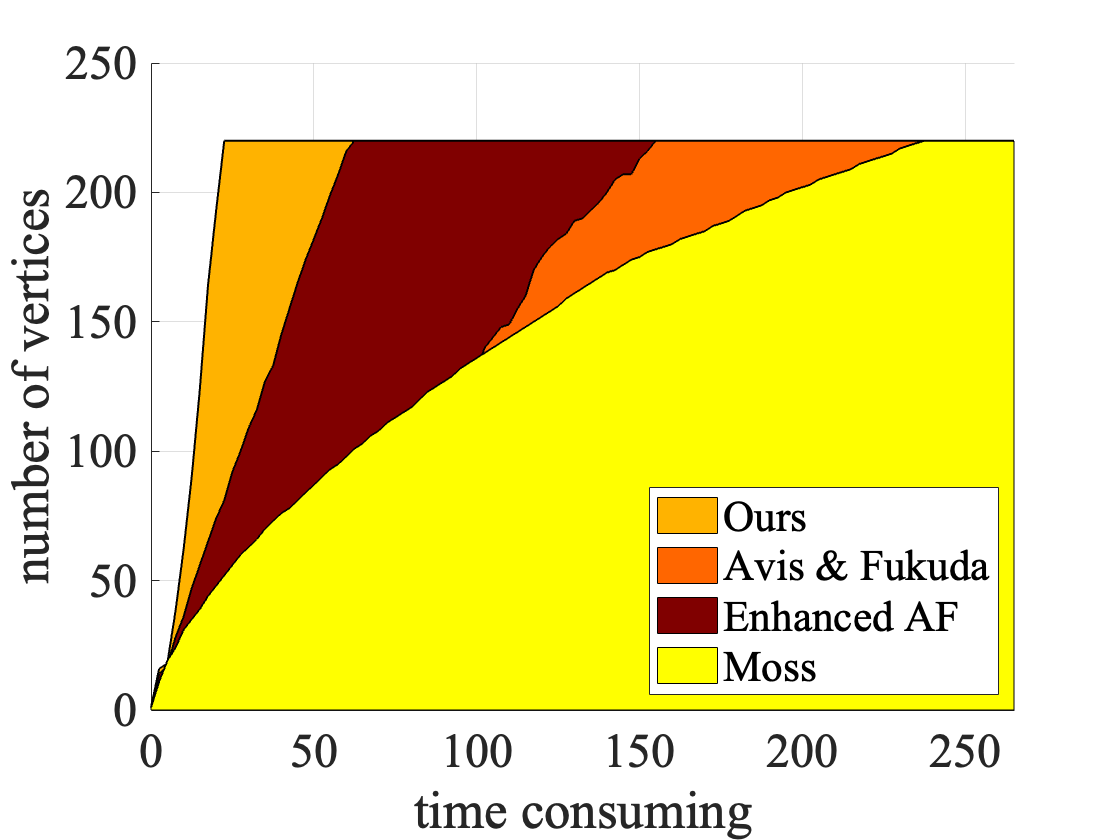}
    \end{minipage}%
    \hspace{0.5cm}
    \begin{minipage}{0.45\linewidth}
    \includegraphics[width=\linewidth]{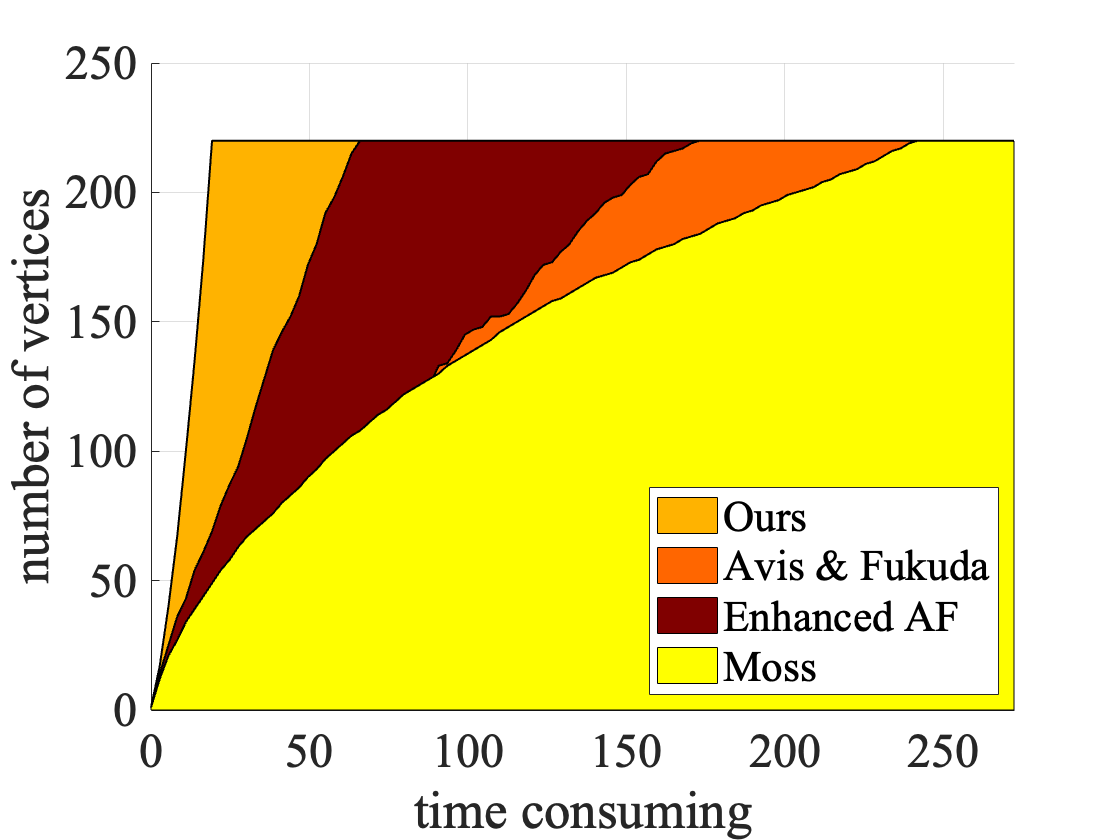}
    \end{minipage}
    \caption{Vertices in 21-st(LHS) and 41-st(RHS) arrangement}
    \label{fig:2141}
    \vspace{-0.4cm}
\end{figure}

{\textbf{Result Analysis}: As Tables \ref{table:cubes_time}-\ref{table:relu-time} shows, in arrangements with a smaller amount of hyperplanes, both ours and the \texttt{Moss} method demonstrate higher efficiency, with the \texttt{Moss} method even achieving the highest efficiency in some cases. However, as the sample size $n$ and the dimensionality $d$ increase, the efficiency of the \texttt{Moss} method declines significantly, while ours consistently maintains high efficiency. Additionally, when comparing the performance of the \texttt{Enhanced AF} and the original \texttt{AF}, which typically rank second and third, it is clear that the \texttt{Enhanced AF} is superior to the original \texttt{AF}. This finding corroborates the advantages of directly employing the if-and-only-if condition to determine valid reverse Criss-Cross pivots in \hyperref[analysis_iff]{Appendix}.}

Moreover, as Figures \ref{fig:rand45}-\ref{fig:2141} show, for each arrangement, in the beginning, both ours and the \texttt{Moss} method exhibit higher rates in finding vertices than the rest two algorithms. As the search continues, the rates of our algorithm and the \texttt{Enhanced AF} remain relatively stable, the rate of the original \texttt{AF} fluctuates, and the rate of the \texttt{Moss} method gradually declines. Due to the stability, the \texttt{Enhanced AF} and ours enjoy the high computational efficiency than the original \texttt{AF}. The gradual decline in the search rate reveals that the \texttt{Moss} method is handicapped in handling complex cases.

The reason for the stability of the \texttt{Enhanced AF} and ours lies in the if-and-only-if conditions for detecting valid reverse pivots, where the computational cost at each position of any dictionary is the same. In contrast, the variability in efficiency of the original \texttt{AF} arises from the varying computational costs at each position, leading to a significant amount of unnecessary calculations. The decreased efficiency of the \texttt{Moss} can be attributed to the increased time required for checking duplicates, as the number of discovered dictionaries grows. Consequently, it performs well with a small sample size but poorly with a larger sample size.

\section{Conclusion}
\label{sec7}
In this paper, we present the Zero rule for the VE problem of hyperplane arrangements. Compared to the classical Criss-Cross rule, the Zero rule enjoys several computationally friendly properties: i) It can get rid of the objective function, which saves the computation. ii) The condition between the Zero rule and its valid reverse is sufficient and necessary, which is more efficient when determining valid reverse pivots. iii) In $\mathbb{R}^d$, the number of pivot steps is at most $d$, which can transform the unpredictable depth into a certain maximum depth, further saving the computation. iv) Its terminal within the entire arrangement is unique, which is convenient for algorithmic design. Because of these properties, the complexity of the VE algorithm using the Zero rule is $\mathcal{O}(n^2d^2(v-v_d)+ndv_d)$. Moreover, it could be as low as $\mathcal{O}(nd^4v)$ for simple arrangements, which greatly improves the state-of-the-art algorithms. Systematic and comprehensive experiments have confirmed not only the efficiency of using the if-and-only-if-condition to determine valid reverse pivot but also the efficiency of the VE algorithm using the Zero rule. In the future, more efforts should be made to further escalate the algorithmic efficiency.

\section{Acknowledgement}
The authors are grateful for Prof. Ning Cai’s suggestions (Hong Kong University of Science and Technology at Guangzhou), discussions with Mr. Xiangru Zhong, a visiting graduate student at the Department of Mathematics, The Chinese University of Hong Kong, and anonymous reviewers’ advice.

\bibliographystyle{plain}
\bibliography{cas-refs}

\appendix

\section*{Appendix. Analysis about the checking of valid reverse pivot}\label{analysis_iff}

{To enhance the efficiency of the algorithm by Avis and Fukuda \cite{avis1991pivoting}, we provide the necessary and sufficient conditions for determining whether the pivot $(s, r)$ on a dictionary is a valid reverse Criss-Cross pivot, as illustrated in Figure \ref{fig:cc_iff_explanation}.} These conditions represent the minimal requirements to confirm a valid reverse Criss-Cross pivot.

\begin{prop}\label{prop:cc_iff}
Let $x_B=\bar{A}x_N$ be an arbitrary dictionary with the objective function. Then pivot $(s,r)$, where $s\in B_{\neq f}$ and $r\in N_{\neq g}$, is a valid reverse Criss-Cross pivot if and only if either (a)
\begin{align*}
&(a1)~\bar{a}_{sg}>0,\bar{a}_{sr}>0,\bar{a}_{sj}\geq0 \ for \ j\in N_{\neq g},j<s.\\
       &(a2)~\forall j<r, \ \mathrm{ If } \ j\in B_{\neq f}, \bar{a}_{jg}\bar{a}_{sr}\geq\bar{a}_{jr}\bar{a}_{sg}; \ \mathrm{ If } \ j\in N_{\neq g},\bar{a}_{fr}\bar{a}_{sj}\geq \bar{a}_{fj}\bar{a}_{sr}.\\
       &(a3)~\textup{ If } s<r, \bar{a}_{fr}\leq0.
\end{align*}
or (b)
\begin{align*}
    &(b1)~\bar{a}_{fr}<0,\bar{a}_{sr}<0,\bar{a}_{ir}\leq0 \ \mathrm{ for }  \ i\in B_{\neq f},i<r.\\
    &(b2)~\forall i<s, \ \mathrm{ If } \ i\in B_{\neq f}, \bar{a}_{ir}\bar{a}_{sg}\geq\bar{a}_{ig}\bar{a}_{sr}; \ \mathrm{ If } \ i\in N_{\neq g},\bar{a}_{fi}\bar{a}_{sr}\geq \bar{a}_{fr}\bar{a}_{si}.\\
    &(b3)~\mathrm{ If }~ r<s, \bar{a}_{sg}\geq0 .
\end{align*}
holds. 
\end{prop}

\noindent\textit{Proof}: Let $\x_{\tilde{B}}=\tilde{A}\x_{\tilde{N}}$ be the dictionary resultant from pivoting $(s,r)$ on $\x_B=\bar{A}\x_N$.

$\Rightarrow$: Assume $(s,r)$ be a valid reverse Criss-Cross pivot. Then the Criss-Cross rule pinpoint $(r,s)$ on $\x_{\tilde{B}}=\tilde{A}\x_{\tilde{N}}$ to pivot. Assume $i$ to be the smallest index such that $x_i$ is dual or primal infeasible, we study cases $i\in \tilde{B}$ and $i\in \tilde{N}$ separately.

\underline{\textit{Case 1}}. Let $i\in \tilde{B}$, then $r=i$ and $\tilde{a}_{rg}<0$, $s$ be the smallest index such that $\tilde{a}_{rs}>0$.

(a1): $\bar{a}_{sr}=\frac{1}{\tilde{a}_{rs}}>0$ and $\bar{a}_{sg}=-\tilde{a}_{rg}\bar{a}_{rs}>0$. Since $\forall j\in \tilde{N}_{\neq g}$ and $j<s$, $\tilde{a}_{rj}\geq0$, then $\bar{a}_{sj}=-\tilde{a}_{rj}\bar{a}_{rs}\geq0$. Note that $\bar{a}_{sr}\geq0$, then $\forall j\in N_{\neq g} $ and $ j<s,\bar{a}_{sj}\geq 0$, (a1) holds.

(a2): $\forall j<r, x_j$ is either primal or dual feasible in $\x_{\tilde{B}}=\tilde{A}\x_{\tilde{N}}$. 

Note that as $j=s$, $\bar{a}_{jg}\bar{a}_{sr}=\bar{a}_{jr}\bar{a}_{sg}$. Also note that for those $j\in \tilde{B}_{\neq f}$ and $j<r$, $x_j$ is primal feasible on $\x_{\tilde{B}}=\tilde{A}\x_{\tilde{N}}$, then $\tilde{a}_{jg}=\bar{a}_{jg}-\frac{\bar{a}_{jr}\bar{a}_{sg}}{\bar{a}_{sr}}\geq0 
$ means $  \bar{a}_{jg}\bar{a}_{sr}\geq\bar{a}_{jr}\bar{a}_{sg}$ $\forall j\in B_{\neq f}$ and $j<r$. Since $\forall j\in \tilde{N}_{\neq g}$ and $j<r$, $x_j$ dual feasible on $\x_{\tilde{B}}=\tilde{A}\x_{\tilde{N}}$, then $\tilde{a}_{fj}=\bar{a}_{fj}-\frac{\bar{a}_{fr}\bar{a}_{sj}}{\bar{a}_{sr}}\leq0$ means $ \bar{a}_{fr}\bar{a}_{sj}\geq \bar{a}_{fj}\bar{a}_{sr}$. (a2) holds.

(a3): (a1-a2) are sufficient for $s>r$. In the case of $s<r$, $x_s$ on $\x_{\tilde{B}}=\tilde{A}\x_{\tilde{N}}$ is dual feasible, $\tilde{a}_{fs}\leq0$, and  $\bar{a}_{fr}=\tilde{a}_{fs}\bar{a}_{sr}\leq0$. (a3) holds.

\underline{\textit{Case 2}}. Let $i\in \tilde{N}$, then $s=i$ and $\tilde{a}_{fs}>0$, $r$ be the smallest index such that $\tilde{a}_{rs}<0$. Similar to the study of the first case, (b1-b3) can be obtained.

$\Leftarrow$: Assume one of the (a) and (b) holds.

\underline{\textit{Case 1}}. Let the part (a) hold. Note that $\forall j\in \tilde{B}_{\neq f}$ and $j<r$, by (a1) and (a2), $\tilde{a}_{jg}=\bar{a}_{jg}-\frac{\bar{a}_{jr}\bar{a}_{sg}}{\bar{a}_{sr}}=\frac{1}{\bar{a}_{sr}}(\bar{a}_{jg}\bar{a}_{sr}-\bar{a}_{jr}\bar{a}_{sg})\geq0, x_j$ on $\x_{\tilde{B}}=\tilde{A}\x_{\tilde{N}}$ is primal feasible. Also $\forall j\in \tilde{N}_{\neq g}$ and $j<r$, if $s\geq r$, by (a2), $\tilde{a}_{fj}=\bar{a}_{fj}-\frac{\bar{a}_{fr}\bar{a}_{sj}}{\bar{a}_{sr}}=-\frac{1}{\bar{a}_{sr}}(\bar{a}_{fr}\bar{a}_{sj}-\bar{a}_{fj}\bar{a}_{sr})\leq0, x_j$ is dual feasible on $\x_{\tilde{B}}=\tilde{A}\x_{\tilde{N}}$. If $s<r$, for the additional case $j=s$, by (a3), $\tilde{a}_{fs}=\frac{\bar{a}_{fr}}{\bar{a}_{sr}}\leq0, x_s$ is dual feasible on $\x_{\tilde{B}}=\tilde{A}\x_{\tilde{N}}$. Hence, for any $j<r$, $x_j$ in $\x_{\tilde{B}}=\tilde{A}\x_{\tilde{N}}$ is dual or primal feasible. 

Since $\tilde{a}_{rg}=-\frac{\bar{a}_{sg}}{\bar{a}_{sr}}<0$, $x_r$ in $\x_{\tilde{B}}=\tilde{A}\tilde{N}$ is primal infeasible, the Criss-Cross rule pinpoint at $r\in\tilde{B}_{\neq f}$. Since $\forall j \in \tilde{N}_{\neq g}$ and $j<s$, by (a1), $\tilde{a}_{rj}=-\frac{\bar{a}_{sj}}{\bar{a}_{sr}}\leq0$ and $\tilde{a}_{rs}=\frac{1}{\bar{a}_{sr}}>0$, thus Criss-Cross rule gives $s\in \tilde{N}$.

Hence, the Criss-Cross rule is fixed $(r,s)$ in $\x_{\tilde{B}}=\tilde{A}\x_{\tilde{N}}$, and $(s,r)$ is a valid reverse Criss-Cross pivot.

\underline{\textit{Case 2}}. Let the part (b) hold. Similar to the first case, applying the Criss-Cross rule on $\x_{\tilde{B}}=\tilde{A}\x_{\tilde{N}}$
the first step fix at $s\in\tilde{N}_{\neq g}$ and second step pinpoint $r$ in the column of $\x_s$, it follows $(s,r)$ is still valid reverse Criss-Cross pivot.

\rightline{$\Box$}\vspace{-0.25cm}

\begin{figure}[htbp]
    \centering
    \includegraphics[width=0.65\linewidth]{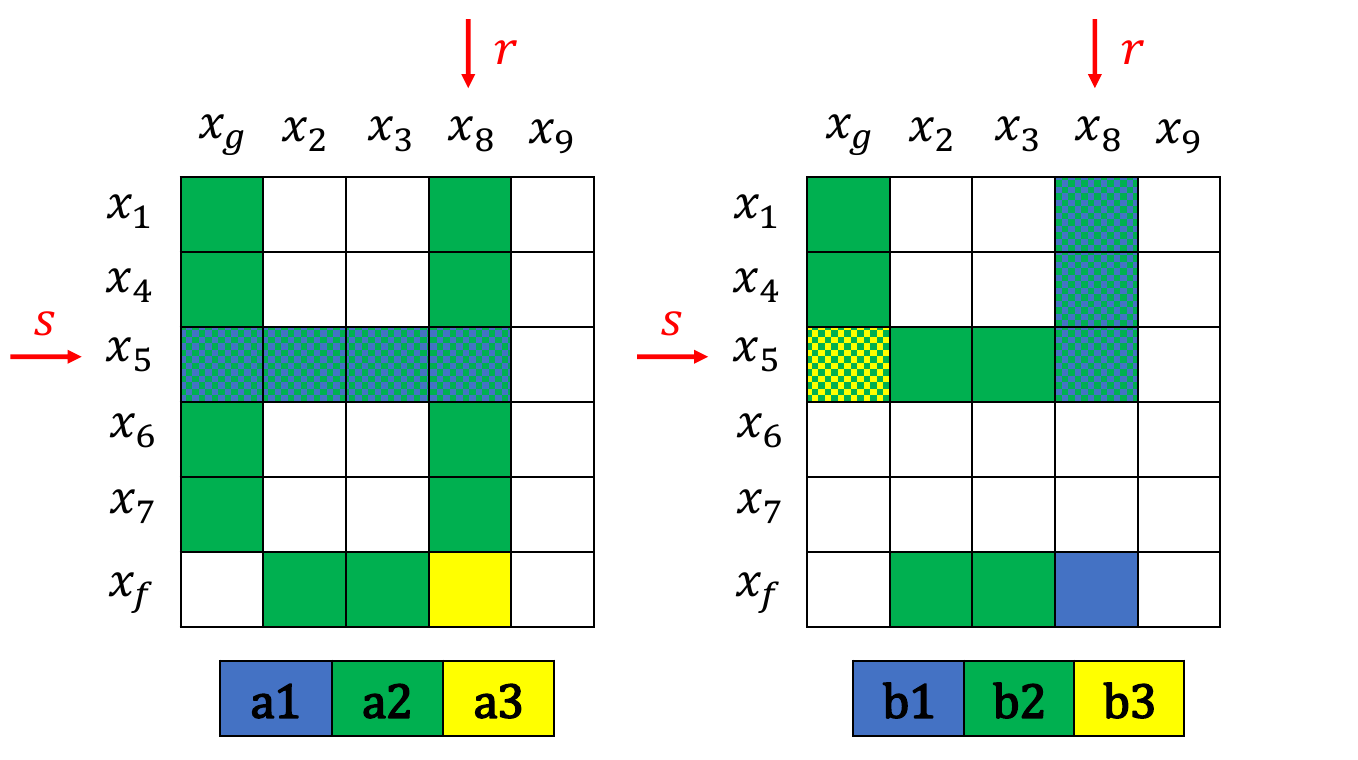}
    \caption{Illustration of the Proposition \ref{prop:cc_iff} to the valid reverse Criss-Cross pivot.}
    \label{fig:cc_iff_explanation}
    \vspace{-0.5cm}
\end{figure}

By applying this condition, checking whether a pivot is a valid reverse Criss-Cross pivot can be simplified as follows: if the necessary and sufficient conditions are met, then it is a valid reverse Criss-Cross pivot. Thus update the dictionary and reset the check position. Conversely, if the conditions are not met, one can proceed to check the next position in the current dictionary. The modification is similar with the RHS of Figure \ref{fig:reverse_check}. The only difference is that the rule checked here is the Criss-Cross rule rather than the Zero rule. A graphic explanation of the \texttt{Enhanced AF} is shown in the LHS of Figure \ref{fig:search_explanation}, and the pseudocode under this formulation is shown in Algorithm \ref{alg:VE_alg_reformulated}. 

Based on the necessary and sufficient conditions and Figure \ref{fig:cc_iff_explanation}, it can be observed that the necessary and sufficient conditions require checking at most $2(2d + 2(n - d)) = 4n$ positions. Thus the analysis for \texttt{Enhanced AF} is relatively straightforward. Regardless of whether the pivots are valid or not, it consumes at most $4n \in \mathcal{O}(n)$ operations at each position in the process similar with the RHS of Figure \ref{fig:reverse_check}, the total operations that the algorithm \texttt{Enhanced AF} requires for
testing valid reverse pivots are $\mathcal{O}(n^2dv)$. {Thus, the total operations cost in checking and performing valid reverse pivot are:}
\begin{equation}
    \begin{aligned}
        \mathcal{O}(n^2dv+ndv)=\mathcal{O}(n^2dv) \label{eq:alg_reformulate_RHS}
    \end{aligned}
\end{equation}

As a result, comparing Eq. \eqref{eq:alg_avis_LHS} with Eq. \eqref{eq:alg_reformulate_RHS}, it can be seen that the total cost of the reverse search on each dictionary by the \texttt{Enhanced AF} is always better than the original \texttt{AF}.




\begin{algorithm}[htbp]
    \caption{\texttt{Enhanced AF}, reverse search of an optimal dictionary}
    \begin{algorithmic} 
        \Function{$\mathtt{Search}(B,N,\bar{A})$}{}\\
            \hspace{0.5cm}$i=1,j=2$;\hfill{\% Use $i, j$ represent the $i$-th row, $j$-th column in $\bar{A}$.}
            \If {$\mathtt{lex-min}(B,N,\bar{A})$==1}\\
                \hspace{1cm}Print $B$;
            \EndIf
            \While{$j\leq \mathtt{length}(N)+1$}
                \If{$j\leq \mathtt{length}(N)$}
                    \If{$\mathtt{reverse}(B,N,\bar{A},i,j)==1$} \hfill{\% A valid reverse pivot is found.} \\
                        \hspace{2.1cm}$[B,N,\bar{A}]=\mathtt{pivot}(B,N,\bar{A},i,j)$;\hfill{\% Compute and update dictionary.} 
                        \If {$\mathtt{lex-min}(B,N,\bar{A})==1$}\\
                            \hspace{2.6cm}Print $B$;
                        \EndIf\\
                        \hspace{2.1cm}$i=1,j=2;$\hfill{\% Initialize the position.}
                    \Else\\
                        \hspace{2.1cm}$[i,j]=\mathtt{increment}(i,j)$;\hfill{\% Invalid, go next position.}
                    \EndIf
                \Else\hfill{\% Each position on this dictionary has been checked.}\\
                    \hspace{1.6cm}$[i,j]=\mathtt{select}(B,N,\bar{A})$;\hfill{\% Use Criss-Cross rule to select a position.}
                    \If{both $i$ and $j$ nonempty}\hfill{\% The Criss-Cross rule find a position.}\\
                        \hspace{2.1cm}$[B,N,\bar{A}]=\mathtt{pivot}(B,N,\bar{A},i,j)$;\hfill{\% Back to parent dictionary.}\\
                        \hspace{2.1cm}$[i,j]=\mathtt{find}(i,j)$;\hfill{\% Find the corresponding position.}\\
                        \hspace{2.1cm}$[i,j]=\mathtt{increment}(i,j)$;
                    \Else\hfill{\% The Criss-Cross rule does not find a position.}\\
                        \hspace{2.1cm}Break;\hfill{\% The function $Search$ terminated.}
                    \EndIf
                \EndIf
            \EndWhile
        \EndFunction \vspace{0.2cm}
        
        \noindent\textbf{function} $\mathtt{lex-min}(B,N,\bar{A})$\hfill{\% Test if the current basis output.}\\
        \noindent\textbf{function} $\mathtt{reverse}(B,N,\bar{A},i,j)$\hfill{\% Test if pivot$(B(i),N(j))$ a valid reverse pivot.}\\
        \noindent\textbf{function} $[B,N,\bar{A}]=\mathtt{pivot}(B,N,\bar{A},i,j)$\hfill{\% Compute the new dictionary.}\\
        \noindent\textbf{function} $[i,j]=\mathtt{increment}(i,j)$\hfill{\% Go next position on the dictionary.} \\
        \noindent\textbf{function} $[i,j]=\mathtt{select}(i,j)$\hfill{\% Use the Criss-Cross rule to find a position.}\\
        \noindent\textbf{function} $[i,j]=\mathtt{find}(i,j)$\hfill{\% Find corresponding position after Criss-Cross pivot.}
    \end{algorithmic}
    \label{alg:VE_alg_reformulated}
\end{algorithm}

\end{document}